\documentclass[12pt, a4paper, onecolumn]{article}
\usepackage{PTSans}
\usepackage[charter,cmintegrals,vvarbb,scaled=1.15]{newtxmath}
\DeclareMathAlphabet{\mathcal}{OMS}{cmsy}{m}{n}
\usepackage[]{algorithm2e}
 
\usepackage{amsmath,bm}
\usepackage{graphicx}
\graphicspath{ {images/} }
\usepackage{natbib}
\usepackage{textcomp}
\usepackage{stmaryrd}
\usepackage{float}
\usepackage{tabularx}
\usepackage[font={sf,bf}, position=t,singlelinecheck=off]{subfig}
\usepackage[labelfont=bf, labelsep=period]{caption}
\usepackage{color}
\usepackage{epstopdf}
\usepackage{upgreek}
\usepackage[nottoc]{tocbibind}

\usepackage[percent]{overpic}
\usepackage{lineno}
\usepackage{pifont}
 \usepackage{relsize}

\usepackage[intoc]{nomencl}
\makenomenclature

\usepackage{ifthen}
\usepackage{etoolbox}
\renewcommand\nomgroup[1]{%
  \item[\bfseries
  \ifstrequal{#1}{G}{Greek symbols}{% 
Track changes is off
Everyone
Track changes for everyone
You
Track changes for You
eirik.keilegavlen
Track changes for eirik.keilegavlen
Guests
Track changes for guests
Current file
Overview
4

  \ifstrequal{#1}{R}{Roman symbols}{%
  \ifstrequal{#1}{S}{Other symbols}{%
  \ifstrequal{#1}{A}{Abbreviations}{%
  \ifstrequal{#1}{O}{Operators, dimensionless groups and special functions}{}}}}}%
]}

\newcommand{\ms}[1]{\textcolor{red}{#1}}

\usepackage{geometry}
\usepackage{setspace}
\usepackage{fancyhdr}

\begin{document}
\noindent\textbf{Modelling and discretization of flow in porous media with thin, full-tensor permeability inclusions}
\\
M. Starnoni$^{1,2}$, I. Berre$^1$, E. Keilegavlen$^1$, \& J.M. Nordbotten$^1$
\\
$^1$Department of Mathematics, University of Bergen, Bergen, Norway \\
$^2$Department of Environment, Land and Infrastructure Engineering, Politecnico di Torino, Torino, Italy 
\
\section*{Abstract}
When modelling fluid flow in fractured reservoirs, it is common to represent the fractures as lower-dimensional inclusions embedded in the host medium. Existing discretizations of flow in %fractured 
porous media with thin inclusions
%(which we will in this paper refer to as "local")
assume that the principal directions of the inclusion permeability tensor are aligned with the inclusion orientation.
While this modelling assumption works well with tensile fractures, 
it may fail in the context of faults, where the damage zone surrounding the main slip surface may introduce anisotropy that is not aligned with the main fault orientation.
In this paper, we introduce a generalized dimensional reduced model
%termed "semi-local" as it local in physical space but non-local in modeling space,
which preserves full-tensor permeability effects also in the out-of-plane direction of the inclusion.
The governing equations of flow for the lower-dimensional objects are obtained through vertical averaging.
We present a framework for discretization of the resulting mixed-dimensional problem, 
aimed at easy adaptation of existing simulation tools.
We give numerical examples that show the failure of existing formulations when applied to anisotropic faulted porous media, and go on to show the convergence of our method in both 2D and 3D.
\\
\\
\noindent\textbf{Key points}
\begin{itemize}
\item Existing local discretizations of flow in fractured porous media fail in modelling out-of plane anisotropic properties of thin inclusions
\item We present a new framework to modelling and discretizing flow in porous media with thin, full-tensor permeability inclusions
\item We show convergence of our method in both 2D and 3D faulted porous media
\end{itemize}

\noindent\textbf{Keywords}
discretization, faults, permeability, mixed-dimensional, flow, porous media

\section{Introduction}
Modeling and simulation of flow in porous media with faults, fractures, and other thin inclusions representing discontinuities is central to a wide range of subsurface engineering applications, including geothermal energy exploitation \citep{bodvarsson1982injection}, shale gas extraction \citep{cao2016fully}, carbon sequestration \citep{johnson2009hydraulic}, and energy storage \citep{nagelhout1997investigating}.

%IB: I switched the McClure and Horne (2014) reference with another one, as the McClure one is more focusing on mechanics than flow. 
%We are interested in subsurface engineering applications concerning flow in three-dimensional porous media with thin inclusions, where the latter can be, for example, faults, fractures,and aquifers.  Fluid flow in reservoirs of such configuration is an important process in sev-eral subsurface applications,  including geothermal energy exploitation (B ̈odvarsson andTsang, 1982), shale gas extraction (Cao et al., 2016), carbon sequestration (Johnson et al.,2009), and energy storage (Nagelhout and Roest, 1997)

The inclusions are characterized by a high aspect ratio, and permeability significantly different from that of the
%d-dimensional %EK: Introduced d-dimensional, not sure if we had notation for dimensions already.
host medium; hence, they severely affect flow patterns. This poses a challenge for traditional simulation models, which are based on upscaling of fine-scale details into an equivalent permeability \citep{oda1985permeability,farmer2002upscaling,liu2016mathematical,saevik20133d}. We instead focus on an alternative approach, which explicitly represents the inclusions in the mathematical and simulation models and thereby to a large degree avoids challenges relating to parameter upscaling. To avoid elongated cells at the inclusion in the computational grid, it is common to represent the inclusions as co-dimension one objects embedded in the host medium \citep{boon2018robust, nordbotten2018unified}.
The intersection of inclusions further gives rise to line and point intersections of co-dimension two and three. %dimensions $d-2$ and $d-3$.
Each of these objects (matrix, inclusions, and intersection points and lines) are represented as independent subdomains separated by interfaces. We refer to this representation of the geometry as mixed-dimensional.

Governing equations for fluid flow in lower-dimensional representation of the inclusion can be derived by integration in the direction orthogonal to the inclusion.
This leads to a decomposition of the governing equations into an in-plane component that represents flow within the inclusion, and an out-of-plane component that couples flow between the inclusion and the host medium. While the in-plane flow has been modeled with both linear and non-linear, as well as both isotropic and non-isotropic flow models \citep{martin2005modeling, reichenberger2006mixed, brenner2017gradient, brenner2018hybrid}, existing models for the coupling term are limited by an assumption on orthogonal flow between inclusion and host. Reduced order models for flow were also developed for aquifers, leading to the same set of equations, see for instance \citet{bear1979hydraulics}, \citet{yortsos1995theoretical}, and \citet{nordbotten2011geological}.
These existing models will be denoted as "local" in the following, meaning that each partial differential equation (PDE) contains only quantities associated with the subdomain where the PDE is defined.  

Local models generally work well when the inclusion is a joint (tensile fracture). However, inclusions with a more complex geological history may have significantly more complex flow properties in the out-of-plane direction. For instance, the damage zone in the vicinity of faults may exhibit shear fractures, slip surfaces, and/or deformation bands, as summarized in \citet{fossen2007deformation}. These features introduce secondary permeability anisotropy in the damage zone as they tend to have preferred orientations, as shown by both field studies \citep{fossen2005fault, johansen2008internal} and core analysis \citep{hesthammer2000spatial}. This leads to preferential flow directions that are neither parallel nor orthogonal to the main plane. This type of flow cannot be represented by existing models that employ dimension reduction. To the Authors' best knowledge, the only attempt to modeling faults and their surrounding damage zones in a mixed-dimensional framework can be found in \cite{fumagalli2019multi}. However, they still apply local formulations to model the damage zones as lower-dimensional objects which are connected on one side to the fault and on the other side to the rock matrix, hence conceptually seeing the whole fault zone as a multilayer object. An alternative approach would be to implement the fault core as a transmissibility multiplier and the damage zone by modifying the grid permeability in the cells adjacent to the model faults, as illustrated in \citet{wilson2020ranking}. In the following, we will consistently refer to the thin inclusions as faults, notwithstanding that all methods presented herein can be applied to models of fractures and other thin inclusions, however, we expect that the methods proposed are of more importance for faults.

The contribution of this paper is two-fold: First, we present a generalized dimensional reduced model that can preserve full-tensor permeability effects also in the out-of-plane direction of the fault. The resulting reduced equations have a form similar to that of traditional models, however the more general coupling structure leads to additional terms both in the in-plane and out-of-plane equations. These terms, as well as our whole novel formulation, will be denoted as "semi-local" in the following, emphasizing the fact that the new PDEs will contain quantities that, while physically in the same location, from a modeling perspective reside outside the subdomain where the PDE is defined, specifically the internal boundary between the subdomain and its higher dimensional neighbor.

Multiple discretization schemes have been proposed for the local dimensionally-reduced models, including methods based on finite volumes \citep{helmig1997multiphase, karimi2003efficient, sandve2012efficient}, mixed finite elements \citep{martin2005modeling, boon2018robust}, virtual elements \citep{fumagalli2019dual} and mimetic methods \citep{formaggia2018analysis}. A comparison study for all these discretizations of flow in fractured media can be found in \citet{flemisch2018benchmarks} and \cite{berre2020verification} for 2D and 3D flow, respectively. However, the additional terms arising in our formulation bring the semi-local model outside the scope of previously proposed discretization methods. 
The second contribution of the paper is therefore the derivation of discretization schemes for semi-local models. We achieve this in two stages: First, based on the unified framework for discretization of mixed-dimensional problems with local interface laws presented in \citet{nordbotten2018unified}, we present conditions under which any standard discretization scheme for fixed-dimensional problems can be extended to mixed-dimensional problems with semi-local interface laws. Second, we present a concrete discretization approach based on finite volume methods.
 
The paper is organized as follows. In Sec. \ref{sec:modelling}, the mathematical model is presented, first for a domain with a single fault, and then for a general faults configuration. Thereafter, in Sec. \ref{sec:discretization}, the unified discretization is formulated. After presenting simulation results in Sec. \ref{sec:numerical}, concluding remarks are given in Sec. \ref{sec:conclusions}.

\section{Flow modelling in faulted porous media}
\label{sec:modelling}
In this section, the mathematical model for flow in faulted porous media is presented, first for a porous domain containing a single fault (Sections \ref{sec:single_inclusion} and \ref{sec:model_reduction}), and then for a general network of faults (Section \ref{sec:generalnetwork}). For the general case, we also provide the weak formulation of the interface problem (Sections \ref{sec:mixeddimensionalformulation}-\ref{sec:variationalformulation}), which will be useful from the perspective of implementation. 

\begin{figure} []
\centering
\begin{overpic}[width=0.75\textwidth]{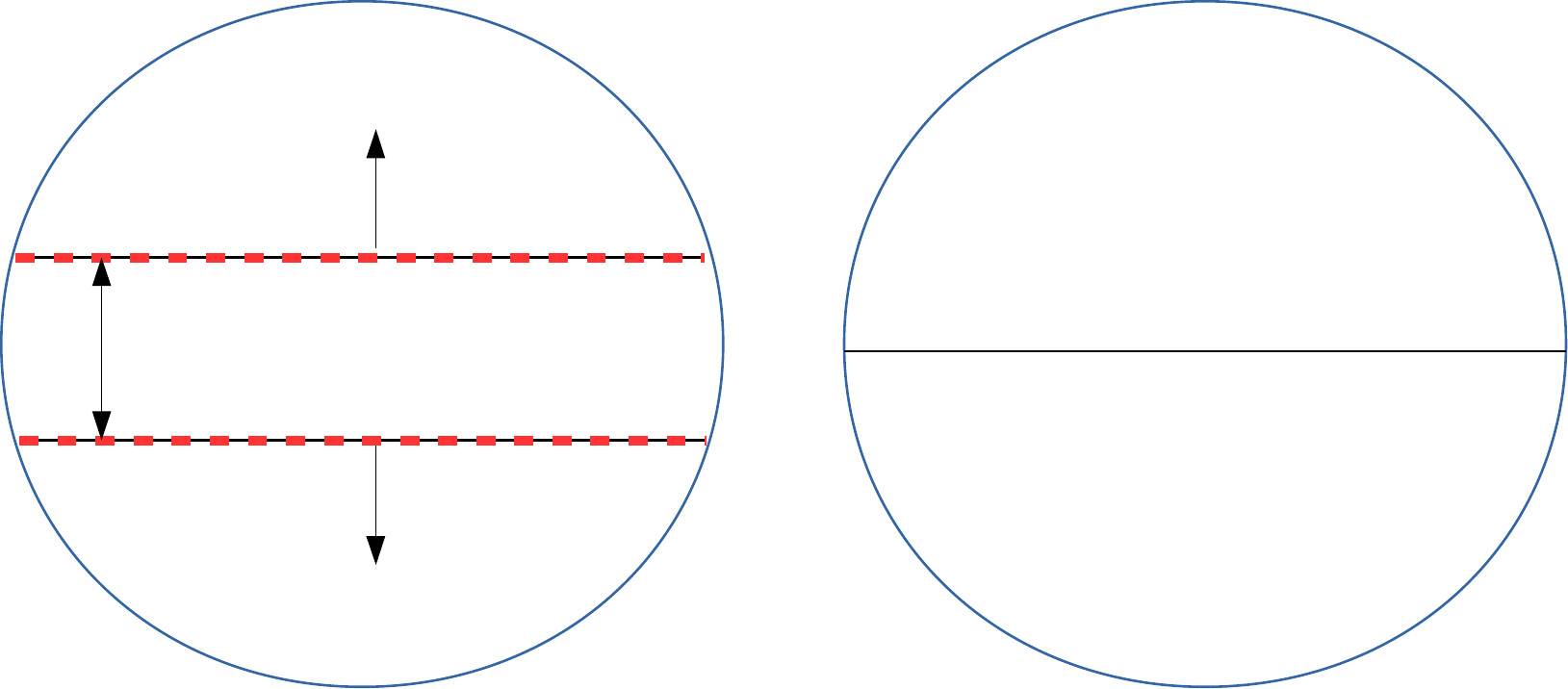}
\put(25,20){\small$\Psi_3$}
\put(30,5){\small$\Psi_2$} 
\put(30,35){\small$\Psi_1$} 
\put(80,23){\small$\Omega_3$}
\put(85,5){\small$\Omega_2$} 
\put(85,35){\small$\Omega_1$} 
\put(23,37){\small${\bm n}_3$} 
\put(23,5){\small${\bm n}_3$} 
\put(10,30){\small$\partial_{\Psi_{3}}\Psi_1$} 
\put(10,12){\small$\partial_{\Psi_{3}}\Psi_2$} 
\put(7,21){\small$a$} 
\end{overpic}
\caption{Representation of the fault as a thin three-dimensional domain $\Psi_3$ (left) and as a two-dimensionl manifold $\Omega_3$ (right). The boundary of $\Psi_j$ adjacent to $\Psi_3$ is denoted by $\partial_{\Psi_3}\Psi_j$, for $j=1,2$, while ${\bm n}_i$ is the normal vector which is always pointing outwards from $\Psi_i$, for $i=1,2,3$.}
\label{fig:upscaling}
\end{figure}
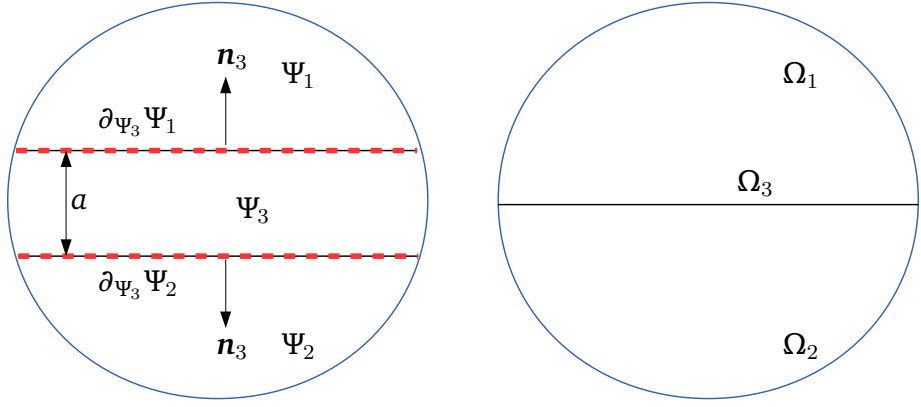

\subsection{Domain with a single fault}
\label{sec:single_inclusion}
We start by considering two three-dimensional porous media $\Psi_1$ and $\Psi_2$, each of them with its Neumann and Dirichlet boundaries $\partial_N$ and $\partial_D$, respectively. The two three-dimensional domains are separated by a fault $\Psi_3$, which is a thin, almost two-dimensional object of thickness $a$ (in the following $a$ will be denoted as the aperture), but which is currently represented as three-dimensional. We note that $\Psi_3$ need not be planar, i.e. $a$ need not be constant. We denote by $\partial_{\Psi_3}\Psi_j$, for $j=1,2$, the boundary of $\Psi_j$ adjacent to $\Psi_3$. Furthermore, let ${\bm n}_i$ be the normal vector which is always pointing outwards from $\Psi_i$. It thus follows that ${\bm n}_3=-{\bm n}_j$ on $\partial_{\Psi_3}\Psi_j$. A representation of the fault as a thin three-dimensional domain $\Psi_3$ is illustrated in the left of Fig. \ref{fig:upscaling}. Darcy flow in the three-dimensional porous medium is then governed by the following set of equations ($i=1,2,3$):
\begin{align}
&\nabla \cdot {\bm q}_i + f_i = 0  \quad &&on \quad \Psi_i \label{eq:mass3}\\
&{\bm q}_i = -{\bm K}_i \nabla p_i \quad &&on \quad \Psi_i \label{eq:darcy3}\\
&\lambda_{3,j} = {\bm q}_3 \cdot {\bm n}_3 = - {\bm q}_j \cdot {\bm n}_j = -\lambda_{j,3} \quad \quad (j=1,2) \quad &&on \quad \partial_{\Psi_{3}}\Psi_j \label{eq:fluxcontinuity}\\
&{\bm q}_i \cdot {\bm n}_i = g_i \quad &&on \quad \partial_N\Psi_i \label{eq:boundaryconditionneumann} \\
&\text{tr }p_i=0 \quad &&on \quad \partial_D\Psi_i \label{eq:3dlast}
\end{align}
Here, $p$ is pressure, $\mathbf{q}$ is the Darcy flux, $f$ is a source, and $\bm K$ is a second-order tensor representing the absolute permeability divided by fluid viscosity. Equation \eqref{eq:mass3} represents mass conservation, while equation \eqref{eq:darcy3} is Darcy's law. Equation \eqref{eq:fluxcontinuity} represents flux continuity conditions on $\partial_{\Psi_3}\Psi_j$, where $\lambda_{3,j}$ represents flow from $\Psi_{3}$ to $\Psi_j$, thus by flux continuity it follows that $\lambda_{3,j}=-\lambda_{j,3}$. Finally, equations \eqref{eq:boundaryconditionneumann}-\eqref{eq:3dlast} are boundary conditions on $\partial_N\Psi_i$ and $\partial_D\Psi_i$, repectively. 

\begin{figure} []
\centering
\includegraphics[width=0.75\textwidth]{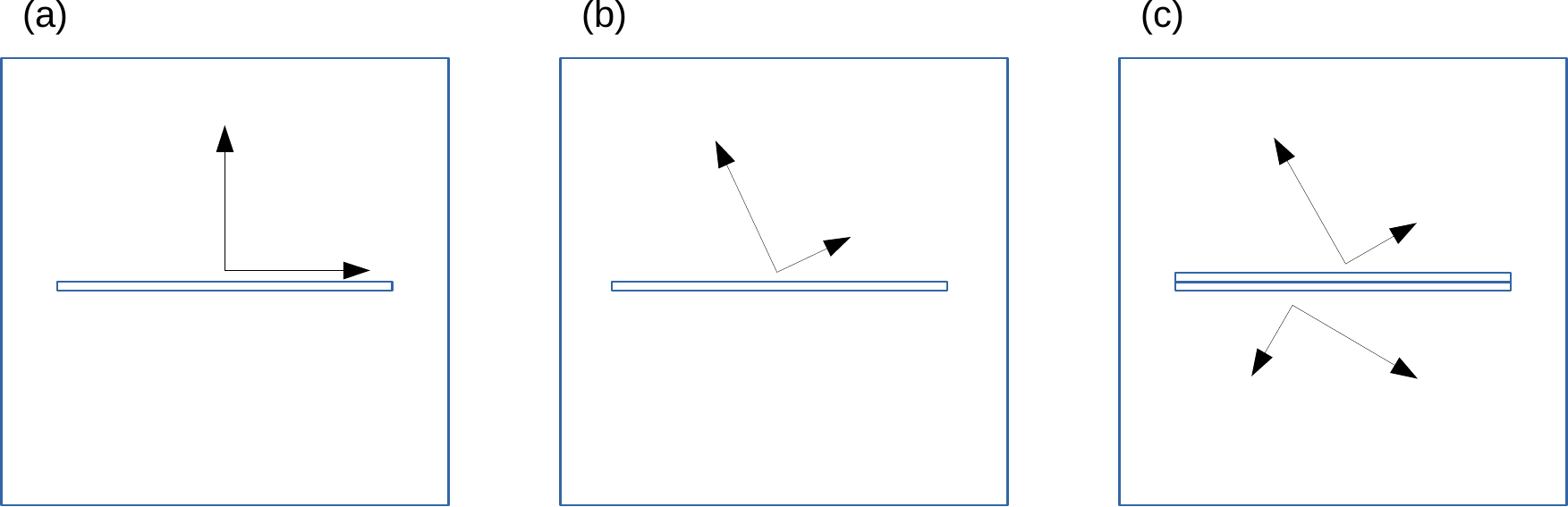}
\caption{Illustration of possible structures of the permeability of a fault embedded in a porous domain indicated by the principal axis of the permeability tensor: (a) orthogonal permeability, (b) homogeneous full-permeability structure, (c) different structure on each half of the fault.}
\label{fig:sketch_permeabilities}
\end{figure}

Before deriving the governing equations for the lower-dimensional manifold, we discuss the decomposition of the permeability tensor within the fault. Existing local laws for faults as embedded thin inclusions assume that the principal directions of the local permeability tensor are aligned with the fault orientation, as illustrated in Fig. \ref{fig:sketch_permeabilities}.a. Hence, more general orientations of the principal permeability directions, shown in Fig. \ref{fig:sketch_permeabilities}.b-\ref{fig:sketch_permeabilities}.c, cannot be represented by existing models. To be concrete, we let the permeability on $\Psi_3$ have the following decomposition in terms of a coordinate system aligned with the fault orientation: 
\begin{equation}
{\bm K}_3=\begin{bmatrix}
{\bm K}_{3,\parallel} & {\bm k}_{3,t}\\
{\bm k}^T_{3,t} & k_{3,\bot}
\end{bmatrix}
\end{equation}
Here, ${\bm K}_{3,\parallel}$ is a $2\times2$ second-order tensor representing the within-fault permeability and $k_{3,\bot}$ is a scalar representing the normal permeability. The off-diagonal term ${\bm k}_{3,t}$ is a two-vector representing the symmetric off-diagonal components of ${\bm K}_3$; for local interface laws, these off-diagonal components are assumed to be negligible, i.e. ${\bm k}_{3,t}=0$ \citep{nordbotten2011geological, berre2020verification}. The inclusion of this anisotropic term leads to significant complications in the modeling and discretization, and is the main topic of this work. With this structure of the fault permeability, the Darcy flux for the fault can be decomposed as ${\bm q}_3 = [{\bm q}_{3,\parallel}, q_{3,\bot}]$, where the 2-vector tangential component ${\bm q}_{3,\parallel}$ and the scalar normal component $q_{3,\bot}$ have the following form:
\begin{align}
\label{eq:darcyTangentialEqui}
&{\bm q}_{3,\parallel}=-{\bm K}_{3,\parallel} \nabla_\parallel p_3 - {\bm k}_{3,t} \nabla_\bot p_3 ,\\
\label{eq:darcyNormalEqui}
&q_{3,\bot} = -{\bm k}_{3,t} \cdot \nabla_\parallel p_3 - k_{3,\bot} \nabla_\bot p_3.
\end{align}
Here, $\nabla_\parallel$ and $\nabla_\bot =\dfrac{\partial}{\partial n}$ represent the in-plane and out-of-plane components of the gradient for the fault, respectively.

\begin{figure} []
\centering
\begin{overpic}[width=0.5\textwidth]{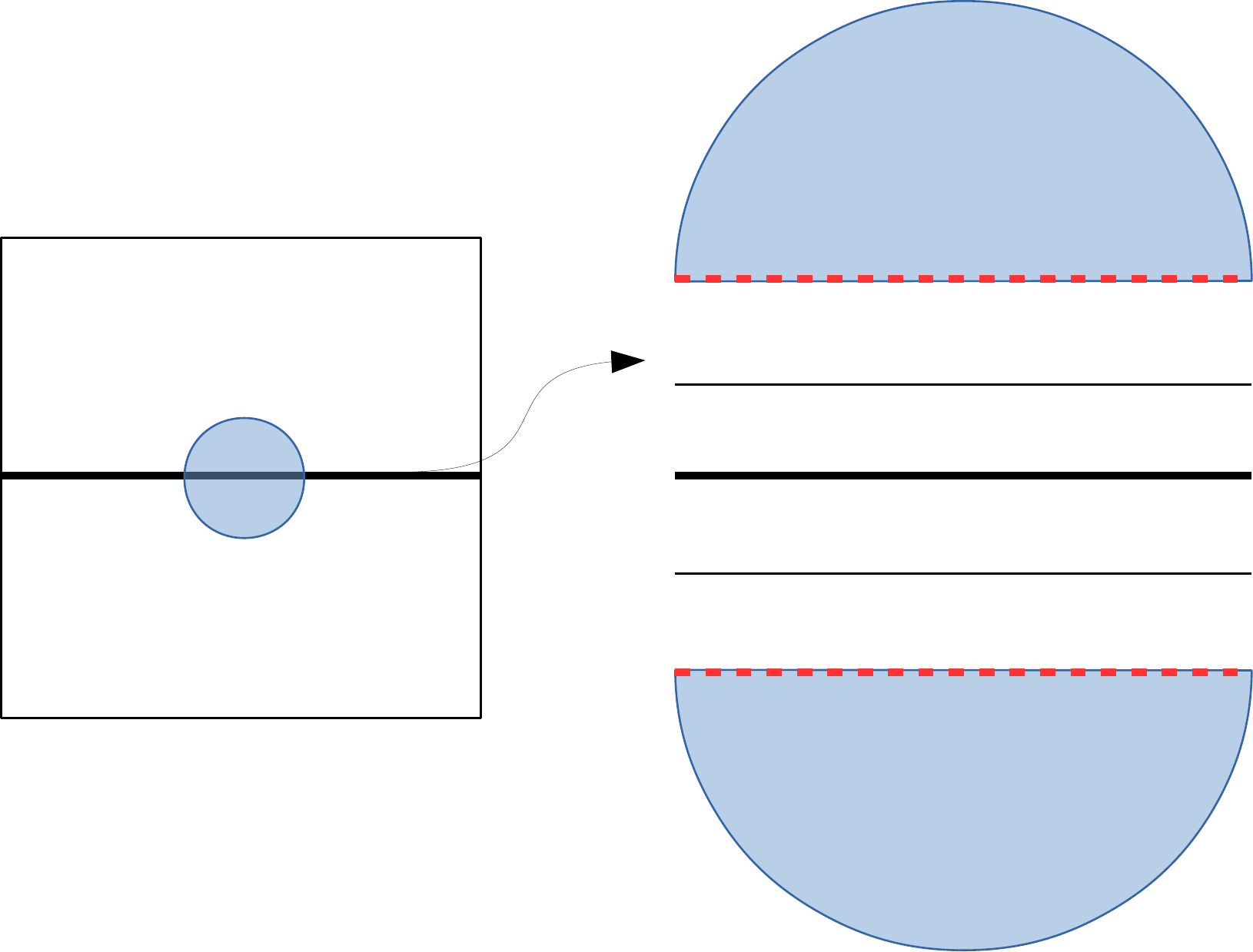}
\put(32,50){\small$\Omega_1$}
\put(32,25){\small$\Omega_2$}
\put(27,40){\small$\Omega_3$} 
\put(102,37){\small$\Omega_3$} 
\put(102,29){\small$\Gamma_{2,3}$} 
\put(102,45){\small$\Gamma_{1,3}$} 
\put(85,56){\small$\partial_{\Omega_3}\Omega_1$} 
\put(85,18){\small$\partial_{\Omega_3}\Omega_2$} 
\put(75,65){\small$\Omega_1$}
\put(75,5){\small$\Omega_2$}
\end{overpic}
\caption{Illustration of the mixed-dimensional geometry. $\Omega_3$ is connected to the higher dimensional neighbors $\Omega_j$ through the interfaces $\Gamma_{j,3}$, for $j=1,2$. Note that $\Omega_3$, $\Gamma_{j,3}$ and $\partial_{\Omega_3}\Omega_j$ are all coinciding in physical space.}
\label{fig:sketch_interfaces}
\end{figure}
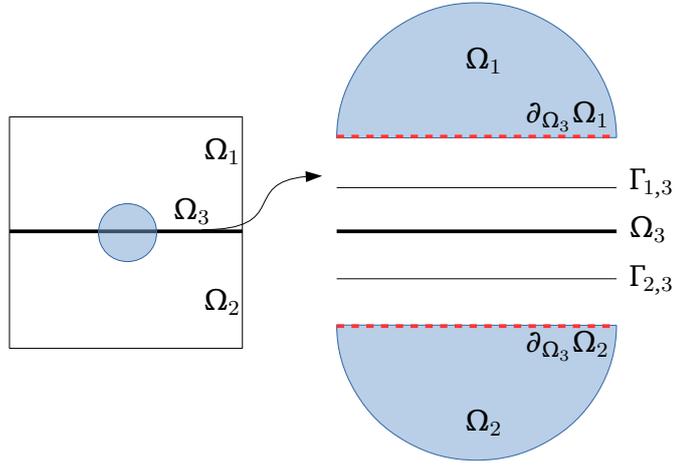

\subsection{Model reduction}
\label{sec:model_reduction}
To proceed, we apply integration over the perpendicular direction to achieve a dimension reduction of the fault. This replaces $\Psi_3$ with a lower-dimensional domain $\Omega_3$ (see right of Fig. \ref{fig:upscaling}). Note that we use $\Psi$ to represent the equi-dimensional geometry, that is all $\Psi_j$ are 3D, and $\Omega$ to denote the mixed-dimensional geometry. We also introduce two interfaces $\Gamma_{j,3}$ on each side $j=1,2$ of $\Omega_3$, as illustrated in Fig. \ref{fig:sketch_interfaces}. The interfaces physically represent the half zone comprised between the fault and either side of the surrounding matrix. In a mixed-dimensional setting, they have no perpendicular extent, and serve as connectors between two objects of different dimensions. Note that, due to the dimension reduction of the model, $\Omega_3$, $\Gamma_{1,3}$, $\Gamma_{2,3}$, $\partial_{\Omega_{3}}\Omega_1$ and $\partial_{\Omega_{3}}\Omega_2$ are all coinciding in physical space. Furthermore, we define the integrated Darcy flux ${\bm q}_3^{(2)}$ and the average pressure $p^{(2)}_3$, respectively as
\begin{equation}
{\bm q}_3^{(2)}=\int_{-a/2}^{a/2}{\bm q}^{(3)}_{3,\parallel} dn, \quad \quad p^{(2)}_3 = \dfrac{1}{a}\int_{-a/2}^{a/2} p_3^{(3)} dn.
\end{equation}
Here, we use subscripts to index the domains, and superscripts (when necessary for clarity) to emphasize the effective topological dimension of the domain, e.g. $p_3^{(3)}$ and $p_3^{(2)}$ are the pressures within the fault in the 3D (on $\Psi_3$) and 2D (on $\Omega_3$) representations, respectively. 
When passing to a mixed-dimensional representation of the geometry, i.e. when integrating eqs. \eqref{eq:mass3} and \eqref{eq:darcyTangentialEqui} along the perpendicular direction, the out-of-plane component of the gradient is converted into a jump operator as follows:
\begin{equation}
\label{eq:jumpOperator}
\int_{-a/2}^{a/2}\nabla_\bot p_3^{(3)} dn = (\text{tr }p_1 - \text{tr }p_2).
\end{equation}
The governing equations for the fault are then obtained from equations \eqref{eq:mass3}, \eqref{eq:darcyTangentialEqui}, \eqref{eq:boundaryconditionneumann} and \eqref{eq:3dlast} by integrating in the perpendicular direction. Moreover, since the fault is assumed to be thin, we assume that the permeability is constant across the perpendicular direction. Together with the definitions above, this results in
\begin{align}
&\nabla_3 \cdot {\bm q}_3^{(2)} - (\lambda_{1,3} + \lambda_{2,3}) + f_3^{(2)}=0 \quad &&on \quad \Omega_3 \label{eq:2dfirst}\\
&{\bm q}_3^{(2)} = - a{\bm K}_{3,\parallel}\nabla_3 p_3^{(2)} + {\bm \mu}_{1,3}+ {\bm \mu}_{2,3} \quad &&on \quad \Omega_3  \label{eq:darcy2D}\\
&{\bm q}_3^{(2)} \cdot {\bm n}_3^{(2)} = g_3^{(2)} \quad &&on \quad \partial_N\Omega_3 \\
&\text{tr }p_3^{(2)}=0 \quad &&on \quad \partial_D\Omega_3 \label{eq:2dlast}
\end{align}
where we have also introduced the integrated source term and boundary flux
\begin{equation}
f_3^{(2)}=\int_{-a/2}^{a/2}f_3^{(3)}dn ,\quad \quad g_3^{(2)}=\int_{-a/2}^{a/2}g_3^{(3)}dn.
\end{equation}
We emphasize that the differential operator $\nabla_3$ in eqs. \eqref{eq:2dfirst}-\eqref{eq:darcy2D} operates on the manifold $\Omega_3$. Compared to traditional upscaled models, see for instance \citet{nordbotten2018unified}, additional terms ${\bm \mu}_{j,3}$ appear in equation \eqref{eq:darcy2D}, analogous to the flux terms $\lambda_{j,3}$ in equation \eqref{eq:2dfirst}. This two-vector term, which is not present in previous work, represents the within-fault flux induced by pressure differences between the fault and the surrounding matrix, and is defined for either side of the fault as
\begin{equation}
\label{eq:jump_definition}
{\bm \mu}_{j,3} = \epsilon_{j,3} {\bm k}_{3,t}(p_3^{(2)}-\text{tr }p_j),
\end{equation}
where the permutation variable $\epsilon_{j,3}$ is positive if the coordinate systems of $\Omega_3$ and $\partial_{\Omega_3}\Omega_j$ coincide, and negative otherwise. 

To complete the model, we derive a constitutive law for $\lambda_{j,3}$. This is obtained by integrating equation \eqref{eq:darcyNormalEqui} in the perpendicular direction, that is
\begin{equation}
\label{eq:darcyNormalEquiIntegrated}
\int_{-a/2}^{a/2}q_{3,\bot}^{(3)} dn=-\int_{-a/2}^{a/2}{\bm k}_{3,t} \cdot \nabla_\parallel p_3^{(3)} dn
-\int_{-a/2}^{a/2}k_{3,\bot} \nabla_\bot p_3^{(3)} dn.
\end{equation}
The left hand side of equation \eqref{eq:darcyNormalEquiIntegrated} is approximated using the trapeizodal rule, that is
\begin{equation}
\label{eq:normalFluxTrapeizodalRule}
\int_{-a/2}^{a/2}q_{3,\bot}^{(3)} dn \approx \dfrac{a}{2}(\epsilon_{1,3}\lambda_{1,3}+\epsilon_{2,3}\lambda_{2,3}),
\end{equation} 
where continuity of the flux across the boundary between the fault and the surrounding matrix is applied. The first term at the right hand side of equation \eqref{eq:darcyNormalEquiIntegrated} is approximated as
\begin{equation}
\label{eq:approximatedGradient}
\int_{-a/2}^{a/2}{\bm k}_{3,t} \cdot \nabla_\parallel p_3^{(3)} dn = {\bm k}_{3,t} \cdot \int_{-a/2}^{a/2} \nabla_\parallel p_3^{(3)} dn \approx \ms{a}{\bm k}_{3,t} \cdot \nabla_3 p_3^{(2)}.
\end{equation}
Finally, the second term at the right hand side of \eqref{eq:darcyNormalEquiIntegrated} is resolved using the jump operator defined in equation \eqref{eq:jumpOperator} as follows:
\begin{equation}
\label{eq:fluxlocalComponent}
\int_{-a/2}^{a/2}k_{3,\bot} \nabla_\bot p_3^{(3)} dn = \epsilon_{1,3} k_{3,\bot} (p_3^{(2)}-\text{tr }p_1) + \epsilon_{2,3} k_{3,\bot} (p_3^{(2)}-\text{tr }p_2).
\end{equation}
By incorporating eqs. \eqref{eq:normalFluxTrapeizodalRule}, \eqref{eq:approximatedGradient} and \eqref{eq:fluxlocalComponent} into equation \eqref{eq:darcyNormalEquiIntegrated}, we identify the flux $\lambda_{j,3}$ having the following form:
\begin{equation}
\lambda_{j,3}= - k_{3,\bot} \dfrac{2(p_3^{(2)}-\text{tr }p_j)}{a}-\epsilon_{j,3}{\bm k}_{3,t} \cdot \nabla_3 p_3^{(2)}. \label{eq:flux}
\end{equation}
Here, the first term on the right-hand side represents the local component of the constitutive law, while the second part is the semi-local contribution that induces a flux across $\Gamma_{j,3}$ due to the pressure gradient within the lower-dimensional manifold $\Omega_3$. 

Inspecting equations \eqref{eq:jump_definition} and \eqref{eq:flux}, we see that both the normal permeability $k_{3,\bot}$ and the off-diagonal permeability ${\bm k}_{3,t}$ are in the reduced model naturally interpreted as properties of the interface $\Gamma_{j,3}$. In the continuation, we will thus generalize the model as derived above, and index these quantities with the interface, i.e. ${\bm k}_{3,t}\rightarrow {\bm k}_{3,j,t}$ and $k_{3,\bot}\rightarrow k_{3,j,\bot}$ are assigned independently to either side of the fault. 

In summary, omitting superscripts for the sake of clarity, we can write the mixed-dimensional equations \eqref{eq:mass3}-\eqref{eq:3dlast}, \eqref{eq:2dfirst}-\eqref{eq:2dlast}, \eqref{eq:jump_definition}, and \eqref{eq:flux} in a unified way, that is for $i=\lbrace 1,2,3 \rbrace$ 
\begin{align}
&\nabla_i \cdot {\bm q}_i - \sum_{j\in{\hat{S}}_i} \lambda_{j,i}  +f_i = 0  \quad &&on \quad \Omega_i \label{eq:massUNi}\\
&{\bm q}_i = -{\bm \kappa}_{i,\parallel} \nabla_i p_i + \sum_{j\in{\hat{S}}_i} \epsilon_{j,i} {\bm \kappa}_{i,j,t}(p_i-\text{tr }p_j)\quad &&on \quad \Omega_i \label{eq:darcyUni}\\
&{\bm q}_i \cdot {\bm n}_i = \lambda_{i,3} \quad \quad (i\neq 3) \quad &&on \quad \partial_{\Omega_3}\Omega_i\\
&\lambda_{j,3}= - \kappa_{3,j,\bot} (p_3-\text{tr }p_j) -{\bm \kappa}_{3,j,t}\cdot \nabla_3 p_3 \quad \quad (j=1,2) \quad &&on \quad \Gamma_{j,3}\\
&{\bm q}_i \cdot {\bm n}_i = g_i \quad &&on \quad \partial_N\Omega_i \\
&\text{tr }p_i=0 \quad &&on \quad \partial_D\Omega_i \label{eq:lastUni}
\end{align}
where ${\hat{S}}_i$ is the set of neighbors of $\Omega_i$ of higher dimension, e.g. ${\hat{S}}_3=\lbrace\Omega_1,\Omega_2\rbrace$. Equations \eqref{eq:massUNi}-\eqref{eq:lastUni} are complemented with the natural convention that there is no four-dimensional domain in the model, thus ${\hat{S}}_i=\emptyset$ for $i=1,2$, and one clearly has for these three-dimensional domains also that $a=1$, ${\bm K}_\parallel=\bm K$ and $\nabla_i=\nabla$. 

We remark that due to the model reduction, the within-fault permeability ${\bm K}_{3,\parallel}$ and the normal permeability $k_{3,j,\bot}$ scale with the aperture $a$ and its inverse, respectively, while the off-diagonal permeability ${\bm k}_{3,j,t}$ remains as in the equi-dimensional model. In order to present equations \eqref{eq:massUNi}-\eqref{eq:lastUni} without reference to this small parameter, these scalings have been incorporated directly into the material constants. Thus, the mixed-dimensional permeability ${\bm \kappa}_3$ is related to the equi-dimensional ${\bm K}_3$ as follows 
\begin{equation}
{\bm \kappa}_3 = \begin{bmatrix}
{\bm \kappa}_{3,\parallel} & {\bm \kappa}_{3,j,t}\\
{\bm \kappa}^T_{3,j,t} & \kappa_{3,j,\bot}
\end{bmatrix}=\begin{bmatrix}
a{\bm K}_{3,\parallel} & {\bm k}_{3,j,t}\\
{\bm k}^T_{3,j,t} & 2a^{-1}k_{3,j,\bot}
\end{bmatrix}.
\end{equation}
We point out that, when one reduces multiple dimensions at once, these scalings get exponents corresponding to the number of dimensions below the ambient dimension. We also emphasize that the normal and off-diagonal permeabilities are in principle not a property of the fault itself, but instead a property which belongs to the internal interface $\Gamma_{j,3}$ between the fault and either side of the higher-dimensional neighbors. This represents an important extension of the existing local laws for fractured porous media, making the model also applicable to faulted porous media, since it allows for capturing the anisotropic character of the fault damage zone. Moreover, since different values of ${\bm k}_{3,j,t}$ and $k_{3,j,\bot}$ can be assigned to each side of the fault, our model can represent different permeability structures on each side of the fault.

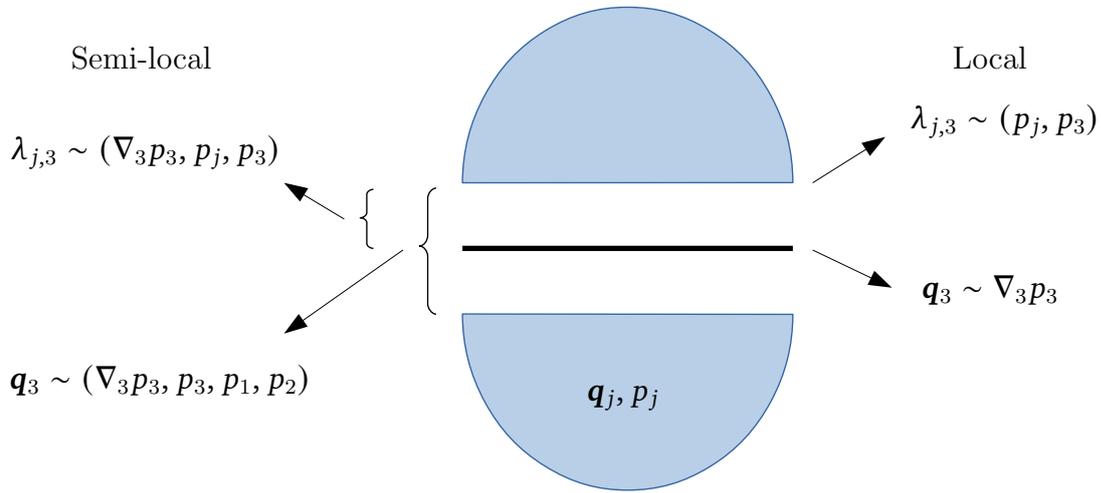
\begin{figure} []
\centering
\begin{overpic}[width=0.5\textwidth]{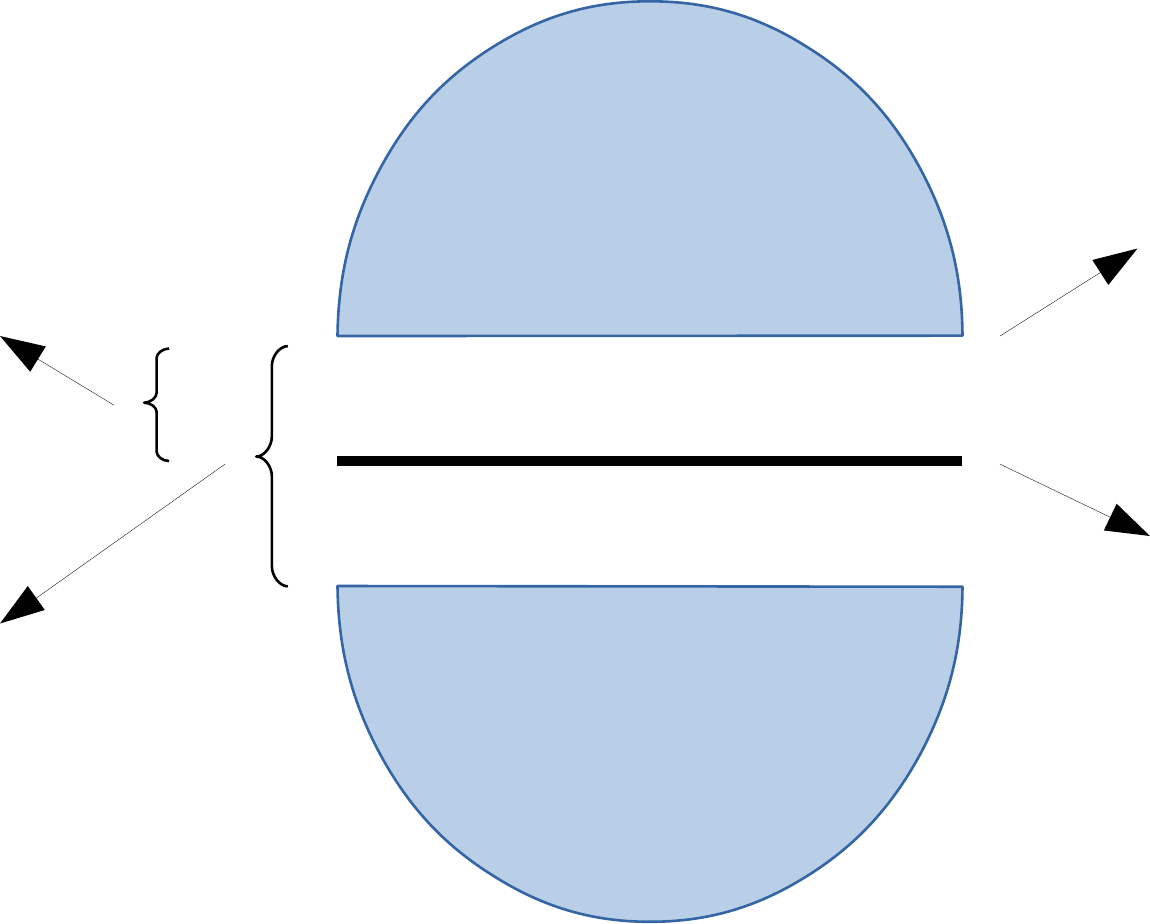}
\put(-35,70){Semi-local}
\put(110,70){Local}
\put(50,15){\small${\bm q}_j,p_j$}
\put(-45,55){\small$\lambda_{j,3}\sim (\nabla_3 p_3, p_j, p_3)$}
\put(-45,17){\small${\bm q}_3 \sim (\nabla_3 p_3, p_3, p_1, p_2)$}
\put(103,60){\small$\lambda_{j,3}\sim (p_j, p_3)$}
\put(105,32){\small${\bm q}_3 \sim \nabla_3 p_3$}
\end{overpic}
\caption{Illustration of the quantities associated with the local and semi-local formulations.}
\label{fig:local_vs_semi_local}
\end{figure}

A schematic illustration of the different quantities and their domain of definition for the local and semi-local formulations is shown in Fig. \ref{fig:local_vs_semi_local}.

\subsection{Domain with a general network of faults}
\label{sec:generalnetwork}
\begin{figure} []
\centering
\begin{overpic}[width=0.4\textwidth]{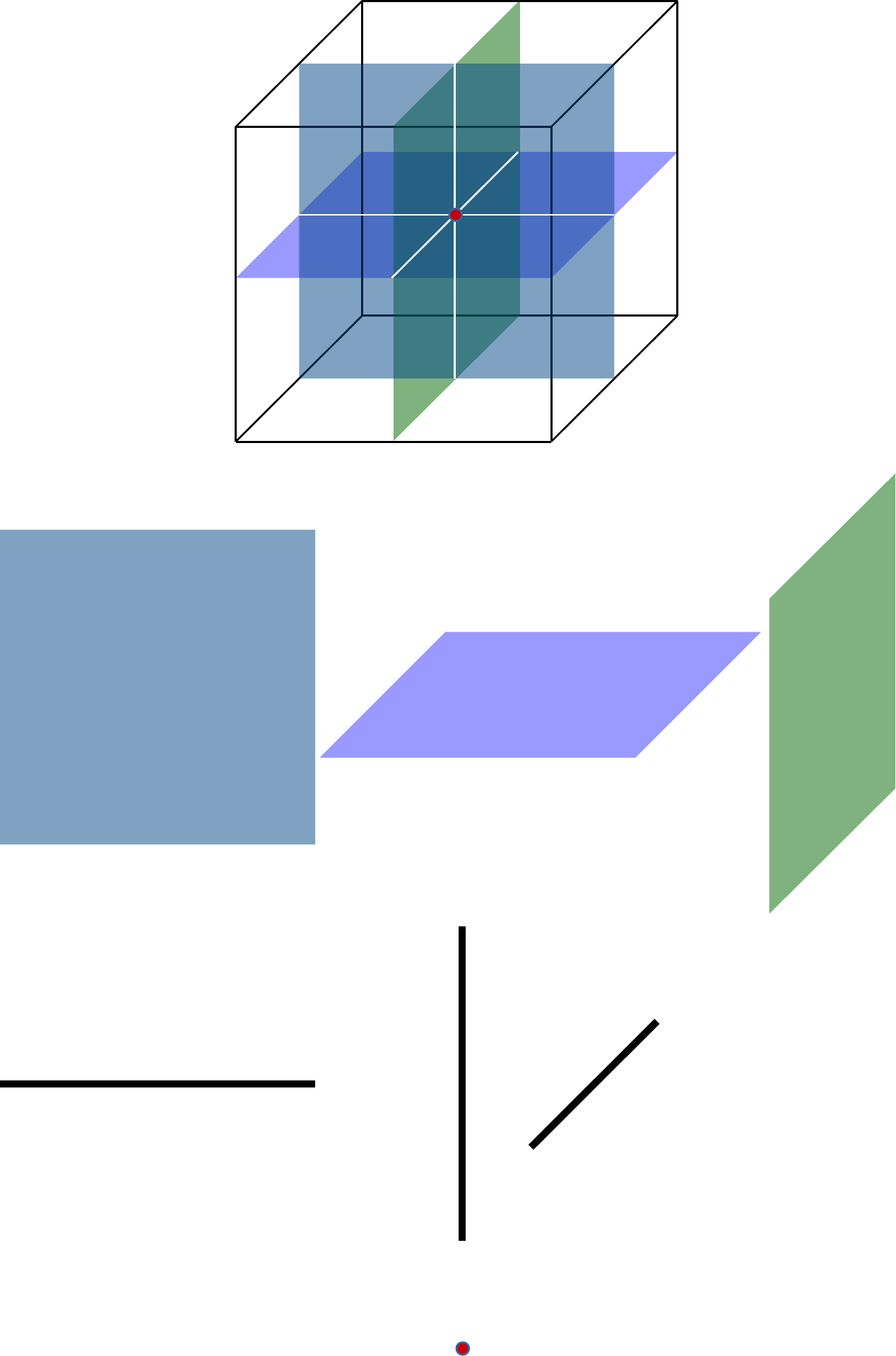}
\put(52,95){\small$\Omega_1$}
\put(5,52){\small$\Omega_2$}
\put(40,47){\small$\Omega_3$} 
\put(60,45){\small$\Omega_4$}
\put(25,20){\small$\Omega_5$}
\put(50,26){\small$\Omega_7$}
\put(36,30){\small$\Omega_6$}
\put(36,2){\small$\Omega_8$}
\end{overpic}
\caption{Illustration of a 3D domain $\Omega_1$ containing three faults $\Omega_j$ ($j=2,3,4$) with their respective three 1D line intersections $\Omega_k$ ($k=5,6,7$) and one 0D point intersection $\Omega_8$.}
\label{fig:sketch_domain_network_inclusions_separate}
\end{figure}
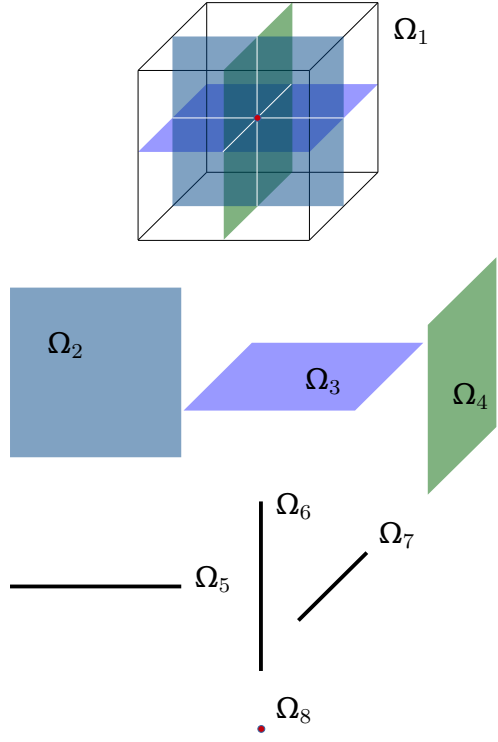

Following the theory by \citet{boon2018robust}, equations \eqref{eq:massUNi}-\eqref{eq:lastUni} can be generalized also to faults intersections, both the one-dimensional (1D) line intersections between two faults and the zero-dimensional (0D) point intersections of three faults (see Fig. \ref{fig:sketch_domain_network_inclusions_separate} for an illustration of the mixed-dimensional geometry). To this end, we use subscripts to index each domain (matrix, fault, or intersection) by number as in the previous section, and let $I$ denote the index set of all domains. 
Superscripts for the topological dimension associated with each individual domain will be consistently omitted, keeping in mind that the dimension is always a property of the domain, i.e. $d=d_i$. Hence, we can write for all $i\in I$ the equations
\begin{align}
&\nabla_i \cdot {\bm q}_i - \sum_{j\in{\hat{S}}_i} \lambda_{j,i}  +f_i = 0  \quad &&on \quad \Omega_i \label{eq:massGeneral}\\
&{\bm q}_i = -{\bm \kappa}_{i,\parallel} \nabla_i p_i + \sum_{j\in{\hat{S}}_i} \epsilon_{j,i} {\bm \kappa}_{i,j,t}(p_i-\text{tr }p_j)\quad &&on \quad \Omega_i \label{eq:darcyGeneral}\\
&{\bm q}_i \cdot {\bm n}_i = \lambda_{i,j} \quad \quad (j \in {\check{S}}_i) \quad &&on \quad \partial_{\Omega_j}\Omega_i \label{eq:fluxGeneral}\\
&\lambda_{j,i}= - \kappa_{i,j,\bot} (p_i-\text{tr }p_j) -\epsilon_{j,i}{\bm \kappa}_{i,j,t}\cdot \nabla_i p_i \quad \quad (j \in {\hat{S}}_i) \quad &&on \quad \Gamma_{j,i} \label{eq:lambdaGeneral}\\
&{\bm q}_i \cdot {\bm n}_i = g_i \quad &&on \quad \partial_N\Omega_i \label{eq:fluxNeumannGeneral}\\
&\text{tr }p_i=0 \quad &&on \quad \partial_D\Omega_i \label{eq:generalLast}
\end{align}
where ${\check{S}}_i$ is the set of neighbors of $\Omega_i$ of lower dimension, e.g. ${\check{S}}_1= \lbrace \Omega_2, \Omega_3, \Omega_4 \rbrace$. It is easy to show that as long as the mixed-dimensional permeabilities are diagonally dominant in the sense of
\begin{equation}
\kappa_{i,j,\bot} \det{\bm \kappa}_{i,\parallel}>{\bm \kappa}_{i,j,t}\cdot{\bm \kappa}_{i,j,t},
\end{equation}
then the coefficients are globally positive definite, and equations \eqref{eq:massGeneral}-\eqref{eq:generalLast} are well-posed as long as $\partial_D\Omega_i$ has non-zero measure for at least one domain \citep{boon2017functional}. 

\subsection{Mixed-dimensional formulation of the fault-matrix flows}
\label{sec:mixeddimensionalformulation}
While equations \eqref{eq:massGeneral}-\eqref{eq:generalLast} constitute a full semi-local model, they are stated in a form which is not immediately amenable for discretization.
This and the following subsection explore the model in more detail, with a goal of rewriting the equations in a form that can be handled by standard discretization schemes with only minimal adaptations. 
A discretization approach based on this reformulation is then given in Section \ref{sec:discretization}. 

In order to simplify the exposition, we will introduce a mixed-dimensional notation following \citet{nordbotten2018unified}. In particular, we will denote the collection of pressure functions as $\mathfrak{p}=\left(p_1,...,p_{|I|}\right)$, and similarly the collection of all fluxes (both in domains and across boundaries) as $\mathfrak{q}=\left({\bm q}_1,...,{\bm q}_{|I|},\lambda_{1,1},...,\lambda_{|I|,|J|}\right)$. It is sometimes convenient to refer explicitly to only the domain or boundary fluxes, and we will therefore sometimes abuse notation and simply write $\mathfrak{q}=(q,\lambda)$. We refer to these as mixed-dimensional functions, and consistently denote them with calligraphic font. We adopt the natural convention that when evaluating a mixed-dimensional function at a point, say $x\in \Omega_i$, then we simply evaluate the function on that domain, so that $\mathfrak{p}(x)=p_i(x)$. In a similar sense, we denote the disjoint union of domains as $\mathfrak{F}=\left(\amalg_i\Omega_i\right) \sqcup \left(\amalg_{j,i}\Gamma_{j,i}\right)$.

With this notion of mixed-dimensional functions, the extension of the divergence and gradient operators to the mixed-dimensional setting is natural. First, we extend the concept of continuous functions by requiring that for $\mathfrak{q}$ to be continuous, then it must hold that, for all $\Gamma_{j,i}$, ${\bm q}_i \cdot {\bm n}_i = \lambda_{i,j}$. Then, for any point $x \in \Omega_i$ we define
\begin{equation}
\left(\mathfrak{D}\cdot \mathfrak{q}\right)(x)=\left[\nabla_i \cdot {\bm q}_i - \sum_{j\in{\hat{S}}_i} \lambda_{j,i} \right]_x
\quad \quad \text{and} \quad \quad
\left(\mathbb{D}\mathfrak{p} \right)(x)=\left[\nabla_ip_i\right]_x,
\end{equation}
while for any point on an interface $x\in\Gamma_{j,i}$ we define
\begin{equation}
\left(\mathbb{D}\mathfrak{p} \right)(x)=\left[p_i-\text{tr }p_j \right]_x.
\end{equation}
Now we can write equations \eqref{eq:massGeneral} - \eqref{eq:generalLast} simply as:
\begin{align}
&\mathfrak{D}\cdot \mathfrak{q}  +\mathfrak{f} = 0  \quad &&on \quad \mathfrak{F} \label{eq:massGeneralGothic}\\
&\mathfrak{q} = -\mathfrak{K}\mathbb{D}p\quad &&on \quad \mathfrak{F} \label{eq:darcyGeneralGothic}\\
&\mathfrak{q} \cdot \mathfrak{n} = \mathfrak{g} \quad &&on \quad \partial_N\mathfrak{F} \label{eq:fluxNeumannGeneralGothic}\\
&\text{tr }\mathfrak{p}=0 \quad &&on \quad \partial_D\mathfrak{F} \label{eq:generalLastGothic}
\end{align}
where we have also introduced the collection of sources $\mathfrak{f}=\left(f_1,...,f_{|I|}\right)$ and the collection of boundary fluxes $\mathfrak{g}=\left(g_1,...,g_{|I|}\right)$. Here, the material coefficients are now all part of the mixed-dimensional permeability $\mathfrak{K}$, which is defined such as that for any mixed-dimensional gradient $\mathfrak{u}=\mathbb{D}p=(u,\mu)$, it holds that for any point $x\in\Omega_i$:
\begin{equation}
\left(\mathfrak{K}\mathfrak{u}\right)(x)={\bm \kappa}_{i,\parallel} u_i - \sum_{j\in{\hat{S}}_i} \epsilon_{j,i} {\bm \kappa}_{i,j,t}\mu_{j,i},
\end{equation}
while for any point on an interface $x\in\Gamma_{j,i}$, it holds that
\begin{equation}
\left(\mathfrak{K}\mathfrak{u}\right)(x)=\kappa_{i,j,\bot} \mu_{j,i} +\epsilon_{j,i}{\bm \kappa}_{i,j,t}\cdot u_i.
\end{equation}
It is then also sometimes convenient to write equation \eqref{eq:darcyGeneralGothic} in matrix form, that is for $\mathfrak{q}=(q,\lambda)$ and $\mathfrak{u}=\mathbb{D}\mathfrak{p}=(u,\mu)$, one has:
\begin{equation}
\label{eq:mixeddimensionaldarcy}
\begin{Bmatrix}
q \\ \lambda
\end{Bmatrix}=-\begin{bmatrix}
\mathfrak{K}_{\Omega\Omega} & \mathfrak{K}_{\Omega\Gamma}\\
\mathfrak{K}_{\Gamma\Omega} & \mathfrak{K}_{\Gamma\Gamma}
\end{bmatrix}\begin{Bmatrix}u \\ \mu\end{Bmatrix}.
\end{equation}

\noindent
Equation \eqref{eq:mixeddimensionaldarcy} highlights the contribution from the semi-local terms in the mixed-dimensional version of Darcy's law.

\subsection{Weak formulation as an interface system}
\label{sec:variationalformulation}
% EK: Rewrote this to something that I understood
The semi-local terms in equations \eqref{eq:massGeneral}-\eqref{eq:generalLast} lead to coupling terms between domains that are local in physical space, but non-local in the mixed-dimensional representation of the geometry.
A critical example are the fault and its sides, which, from the perspective of implementation, we would prefer to only interact via the interfaces $\Gamma_{j,i}$, and not directly, as is the case for the last term in equation \eqref{eq:darcyGeneral}.

Thus we are motivated to consider a reformulation of the governing equations before considering numerical discretizations. We proceed by first performing an LU decomposition of equation \eqref{eq:mixeddimensionaldarcy} as follows:
\begin{equation}
\label{eq:darcyLU}
\mathfrak{K}_U\begin{Bmatrix}q\\ \lambda \end{Bmatrix}
= -\mathfrak{K}_L\begin{Bmatrix}u\\ \mu \end{Bmatrix},
\end{equation}
where $\mathfrak{K}_U$ and $\mathfrak{K}_L$ are defined, respectively, as:
\begin{equation}
\mathfrak{K}_U=\begin{bmatrix}
I & \mathfrak{K}_{\Omega\Gamma}\mathfrak{K}_{\Gamma\Gamma}^{-1}\\ 0 & I
\end{bmatrix} 
\quad \quad \text{and} \quad \quad
\mathfrak{K}_L=\begin{bmatrix}
A_\Omega & 0 \\ \mathfrak{K}_{\Gamma\Omega} & \mathfrak{K}_{\Gamma\Gamma}
\end{bmatrix},
\end{equation}
and $A_\Omega$ is the Schur-complement defined as
\begin{equation}
A_\Omega=\mathfrak{K}_{\Omega\Omega}-\mathfrak{K}_{\Omega\Gamma}\mathfrak{K}_{\Gamma\Gamma}^{-1}\mathfrak{K}_{\Gamma\Omega}.
\end{equation}
Note that, since $\mathfrak{K}_{\Gamma\Gamma}$ consists only of scalar values $(\kappa_{i,j,\bot})$, this reformulation only depends on the trivial inversion of scalars.

In the following it will be helpful to discuss the components of the mixed-dimensional gradient and divergence, and we therefore additionally define the "full jump" $\mathbb{d}\mathfrak{q}$ such that for any point $x\in\Omega_i$ it holds that
\begin{equation}
\label{eq:totaljumpgothic}
\left(\mathbb{d}\mathfrak{q}\right)(x)=\left[-\sum_{j\in{\hat{S}}_i} \lambda_{j,i}\right]_x,
\end{equation}
while the "half jump" $\mathbb{d}^\star \mathfrak{p}$ is simply the restriction of $\mathbb{D}\mathfrak{p}$ to $\Gamma_{j,i}$. We then write (with the natural extension of $\nabla$ and $\nabla\cdot$):
\begin{equation}
\label{eq:pressuregradientGothic}
\mathfrak{D}\cdot\mathfrak{q}=\nabla\cdot q + \mathbb{d}\lambda
\quad \quad \text{and} \quad \quad
\mathbb{D}\mathfrak{p}=\left(\nabla p,\mathbb{d}^\star \mathfrak{p}\right).
\end{equation}
We now proceed by (formally) eliminating internal domain variables, in order to obtain a problem only posed on interfaces. We note that equations \eqref{eq:massGeneralGothic} and \eqref{eq:darcyGeneralGothic} can now be written as the first order system:
\begin{align}
&\mathfrak{D}\cdot\mathfrak{q}=\mathfrak{f} \label{eq:globalmassconservation}\\
&\mathfrak{K}_U \mathfrak{q}=-\mathfrak{K}_L\mathfrak{D}\mathfrak{p} \label{eq:globaldarcy}
\end{align}
where use of equation \eqref{eq:darcyLU} has been made. By writing out equation \eqref{eq:globalmassconservation} in local notation for each $\Omega_i$ and by stating equation \eqref{eq:globaldarcy} explicitly as two equations, we obtain the following set of equations:
\begin{align}
&\nabla\cdot q=f-\mathbb{d}\lambda \label{eq:localmassconservation}\\
&q+A_\Omega\nabla p=-\mathfrak{K}_{\Omega\Gamma}\mathfrak{K}_{\Gamma\Gamma}^{-1}\lambda \label{eq:localdarcy}\\
&\lambda=-\left(\mathfrak{K}_{\Gamma\Omega}\nabla p+\mathfrak{K}_{\Gamma\Gamma}\mathbb{d}^\star \mathfrak{p}\right) \label{eq:localinterfaceflux}
\end{align} 
This reveals that equations \eqref{eq:localmassconservation} and \eqref{eq:localdarcy} form a locally well-posed system (of standard Darcy type) on each $\Omega_i$, and we can therefore consider $p=p(\lambda)$ for any given $\lambda$.

We formalize this concept by introducing the (continuous) solution operators for the standard elliptic value problem on $\Omega_i$, $\mathcal{S}_{\Omega_i}^K$, defined as:
\begin{equation}
\left(\upsilon, \nabla \upsilon,\text{tr }\upsilon,F\right)=\mathcal{S}_{\Omega_i}^K \left(f, \chi, b, \upsilon_0\right),
\end{equation}
where $\upsilon$ is the solution to
\begin{align}
&\nabla\cdot\varphi=f-F \quad \quad &&\text{on} \quad \Omega_i\\
&\varphi=-K\left(\nabla\upsilon + \chi\right) \quad \quad &&\text{on} \quad \Omega_i\\
&\varphi \cdot n =b \quad \quad &&\text{on} \quad \partial\Omega_i \setminus \partial\Omega\\
&\upsilon = 0 \quad \quad &&\text{on} \quad \partial\Omega_i \cap \partial\Omega \\
&\dfrac{1}{|\Omega_i|}\int_{\Omega_i}\upsilon=\upsilon_0 \quad \quad &&\text{if} \quad \partial\Omega_i \cap \partial\Omega \neq \oslash
\end{align}
where $\partial\Omega$ is the global boundary and $F=\dfrac{1}{|\Omega_i|}\left(\int_{\Omega_i}f-\int_{\Omega_i}b\right)$ if $\partial\Omega_i \cap \partial\Omega \neq \oslash$, and zero otherwise. Using this solution operator, we see that the solution to equations \eqref{eq:localmassconservation} and \eqref{eq:localdarcy} can be stated as functions of $\lambda$ (and a set of number of numbers $p_0$ corresponding to the domains where $\partial\Omega_i \cap \partial\Omega \neq \oslash$) as:
\begin{equation}
\label{eq:solutionoperator}
\left(p, \nabla p,\text{tr }p,F\right)_{\Omega_i}\left(\lambda,p_0\right)=\mathcal{S}_{\Omega_i}^{A_i} \left(f_i-\left(\mathbb{d}\lambda\right)_i, A_i^{-1}\left(\mathfrak{K}_{\Omega\Gamma}\mathfrak{K}_{\Gamma\Gamma}^{-1}\lambda\right)_i, \lambda_{\check{I}_i}, p_0\right).
\end{equation}
Inserting $p=p(\lambda,p_0)$ etc. into equation \eqref{eq:localinterfaceflux}, we have now reformulated the fault-matrix problem into a pure interface problem. From the perspective of implementation, we desire to consider the interface problem in the weak sense, and we therefore multiply by test functions $w$ and integrate to obtain the problem: Find $\lambda \in L^2(\Gamma)$ such that, for all $w\in L^2(\Gamma_j)$
\begin{equation}
\label{eq:weakformulation}
\left(\mathfrak{K}_{\Gamma\Gamma}^{-1}\lambda,w\right)_{\Gamma_{j,i}}+\left(\mathfrak{K}_{\Gamma\Gamma}^{-1}\mathfrak{K}_{\Gamma\Omega}\nabla p (\lambda,p_0),w\right)_{\Gamma_{j,i}}+\left(\mathbb{d}^\star \mathfrak{p}(\lambda,p_0),w\right)_{\Gamma_{j,i}}=0
\end{equation}
and $F_i(\lambda,p_0)=0$ if $\partial\Omega_i \cap \partial\Omega \neq \oslash$. We point out that the inner products in equation \eqref{eq:weakformulation} are bounded from a formal perspective, since for $\lambda\in L^2(\Gamma)$, then $p_i\in H^1(\Omega_i)$, and both $\mathfrak{K}_{\Gamma\Gamma}^{-1}\mathfrak{K}_{\Gamma\Omega}\nabla p$ and $\text{tr }p$ will lie in (at least) $L^2(\Gamma_{j,i})$. 

Finally, we emphasize that equations \eqref{eq:solutionoperator}-\eqref{eq:weakformulation} are attractive from the perspective of implementation, since the inner products appearing are easy to evaluate, and the solution operators $\mathcal{S}_{\Omega_i}^{A_i}$ can be approximated by any standard method, as we will detail in the next section.

\section{Discretizations of flow for faulted porous media}
\label{sec:discretization}
% EK: Tried to create a pick-up point for readers who got lost in section 2.4-5.
The equations derived in Section \ref{sec:variationalformulation}, and in particular the interface problem of equation \eqref{eq:solutionoperator}, form the starting point for the discretization approach laid out in this section. 
We present the general discretization framework in Section \ref{sec:unified_discretization}, and discuss implementational aspects in Section \ref{sec:implementation}.

\subsection{Unified discretization}
\label{sec:unified_discretization}
Equation \eqref{eq:solutionoperator} provides a solution operator for the arbitrary standard method used to solve the elliptic boundary value problem \eqref{eq:localmassconservation}-\eqref{eq:localdarcy} on $\Omega_i$. 
To be concrete, we consider each domain $\Omega_i$ and its Neumann boundary $\partial\Omega_i=\partial_N\Omega_i \cup \sum_{j\in\check{S}_i}\partial_{\Omega_j}\Omega_i$ as endowed with a numerical discretization. Then, the solution operator $\mathcal{S}_i$ can be stated as 
\begin{equation}
\mathcal{S}_i:\left[N(\Omega_i),N^{d_i}(\Omega_i),N(\partial\Omega_i)]\rightarrow[N(\Omega_i),N^{d_i}(\Omega_i),N(\partial \Omega_i)\right],
\end{equation}
where $N(\Omega_i)$, $N^{d_i}(\Omega_i)$, and $N(\partial \Omega_i)$ are the discrete representations of $L^2(\Omega_i)$, $\left( L^2(\Omega_i) \right)^{d_i}$, and $L^2(\partial \Omega_i)$, respectively, and ${d_i}$ is the topological dimension of $\Omega_i$. In particular, $\mathcal{S}_i$ takes as input sinks, vector sources, and Neumann data and returns as output pressures, pressure gradients, and pressure traces.
Most discretization schemes for elliptic equations can provide such a solution operator; we discuss the concrete implementation in the next subsection.

To discretize the flux coupling term $\lambda_{j,i}$, we introduce a mortar-like grid $\mathcal{T}_{j,i}$ on the interface $\Gamma_{j,i}$ on which the boundary flux $\lambda_{j,i}$ will be defined. The flux variables are represented as piecewise constant on the mortar grid $\mathcal{T}_{j,i}$, thus $\lambda_{j,i} \in P_0(\mathcal{T}_{j,i})\subset \L^2(\Omega_i)$. In order to allow communications between subdomains, and thus explicitly relate the degrees of freedom of the numerical methods $\mathcal{S}_i$ and the mortar grids $\mathcal{T}_{j,i}$, we introduce projection operators, namely $\Pi_{N(\Omega_i)}$ and $\Pi_{L^2(\Omega_i)}$. The former is the compound operator projecting from the coupling variables on the mortar grids to the subdomain degrees of freedom, that is
\begin{equation}
\begin{split}
\Pi_{N(\Omega_i)}:&\left[L^2(\Omega_i),\left(L^2(\Omega_i) \right)^{d_i},L^2\left(\Omega_{\check{S}_i}\right),L^2(\partial\Omega_i)\right]\\
&\rightarrow\left[N(\Omega_i),N^{d_i}(\Omega_i),N(\partial \Omega_i)\right],
\end{split}
\end{equation} 
while the latter conversely moves from the numerical variables to the coupling variables, that is
\begin{equation}
\begin{split}
\Pi_{L^2(\Omega_i)}:&\left[N(\Omega_i),N^{d_i}(\Omega_i),N(\partial \Omega_i)\right]\\
&\rightarrow\left[L^2(\Omega_i),\left(L^2(\Omega_i) \right)^{d_i},L^2\left(\Omega_{\check{S}_i}\right),L^2(\partial\Omega_i)\right].
\end{split}
\end{equation}
Now, following the variational formulation derived in Sec. \ref{sec:variationalformulation}, we exploit equation \eqref{eq:weakformulation} in order to provide discretization-independent framework for faulted porous media. This takes the form: for given numerical discretizations $\mathcal{S}_i$, find $\lambda_{j,i}\in P_0(\mathcal{T}_{j,i})$, for all $i\in I$ and $j\in\hat{S}_i$ such that
\begin{equation}
\label{eq:darcyNormal}
\begin{split}
&\left(\mathbb{d}^\star \mathfrak{p},w_j\right)_{\Gamma_{j,i}}+\left( \mathfrak{K}^{-1}_{\Gamma\Gamma}\left(\lambda_{j,i}+\mathfrak{K}_{\Gamma\Omega} \cdot \nabla p\right),w_j\right)_{\Gamma_{j,i}}=0 \\
&\quad \text{for all } w_j \in P_0(\mathcal{T}_{j,i})\\
\end{split}
\end{equation}
subject to discrete constraints (for all $i\in I$):
\begin{align}
[p_i,u_i,t_{j}]&=\Pi_{L^2\left(\Omega_i\right)}\mathcal{S}_i(\psi_i + a_i,b_i, c_i) \label{eq:constraintEqu}\\
[a_i,b_i, c_i]&=\Pi_{N(\Omega_i)}\left[-\sum_{j\in\hat{S}_i}\lambda_{j,i},
-\sum_{j\in\hat{S}_i} A_{i}^{-1}\mathfrak{K}_{\Omega\Gamma}\mathfrak{K}_{\Gamma\Gamma}^{-1}\lambda_{j,i}, \sum_{j\in\check{S}_i}\lambda_{i,j}\right] \label{eq:constraintMortar}
\end{align}
where the dummy variables $a_i$, $b_i$ and $c_i$ have the interpretations of sources, forces, and fluxes due to interactions with other domains, respectively. In contrast, the variables $p_i$, $u_i$, and $t_j$ are the pressures, pressure gradients, and pressure traces after projection onto the grids $\mathcal{T}_{j,i}$.

The interpretation of this scheme is as follows. Eq. \eqref{eq:constraintEqu} resolves the internal differential equations in each subdomain, eq.\eqref{eq:constraintMortar} is the projection of variables from the flux grids to the numerical boundary (and source) data, while equation \eqref{eq:darcyNormal} simply states that the flux $\lambda_{j,i}$ between the fault and its surrodundings should satisfy Darcy's law. In the following section, we present the strategy for implementation of this approach and give details for a specific numerical scheme.

\subsection{MPFA discretization}
\label{sec:implementation}
% EK: Refocused this section to first consider general requirements for numerical methods, and next the realization within PorePy.
It is of interest to consider the requirements put on the subdomain solution operators $\mathcal{S}_i$ in some more detail. 
From the variational formulations stated above, we see that for a discretization on a generic subdomain $\Omega_i$ to interact with the interface $\Gamma_j$, we need to provide operators which:
\begin{enumerate}
\item Handle Neumann boundary data on the form $\Pi_{N(\Omega_i)}\lambda_{j}$, for all interfaces $\Gamma_j$ where $\Omega_i$ is the higher-dimensional neighbor.
\item Handle source terms $\Pi_{N(\Omega_i)}\lambda_{j}$ from interfaces $\Gamma_j$ where $\Omega_i$ is the lower-dimensional neighbor.
\item Provide a discrete operator $\text{tr }p_i$ so that $\Pi_{L^2(\Omega_i)}$ can project the pressure trace from $\partial_j \Omega_i$ to $\Gamma_j$ where $\Omega_i$ is the higher-dimensional neighbor.
\item Provide a pressure $p_i$ so that $\Pi_{L^2(\Omega_i)}$ can project the pressure to all $\Gamma_j$ where $\Omega_i$ is the lower-dimensional neighbor.
\item Handle the divergence of vector source terms $\Pi_{N(\Omega_i)}(\nabla \cdot {\bm \mu}_{j,i})$ from interfaces $\Gamma_j$ where $\Omega_i$ is the lower-dimensional neighbor.
\item Provide a pressure gradient $u_i$ so that $\Pi_{L^2(\Omega_i)}$ can project the pressure gradient to all $\Gamma_j$ where $\Omega_i$ is the lower-dimensional neighbor. 
\end{enumerate}

\noindent
The four first requirements are readily available for any discretization scheme for elliptic equations.
Specifically, we have based our solution operators on a cell-centered finite volume method termed the multi-point flux approximation (MPFA) \citep{aavatsmark2002introduction, nordbotten2020introduction}.
Treatment of vector source terms (item 5) is not as natural in primal discretization schemes such as finite elements, but is easy to include in most flux-based discretization methods such as e.g. mixed finite elements. We have employed the approach introduced in \citet{starnoni2019consistent}, which treats the vector source term as part of the discrete divergence operator, and thereby provides an expression of the fluxes in terms of jumps in cell-centers vector sources.
Finally, the pressure gradients are discretized as piece wise constant on each cell from an interpolation of the face cells fluxes (item 6).
We implemented our model in PorePy, an open-source software for simulation of multiphysics processes in fractured porous media  \citep{keilegavlen2020porepy}.

\begin{figure} []
\centering
\begin{overpic}[width=0.25\textwidth]{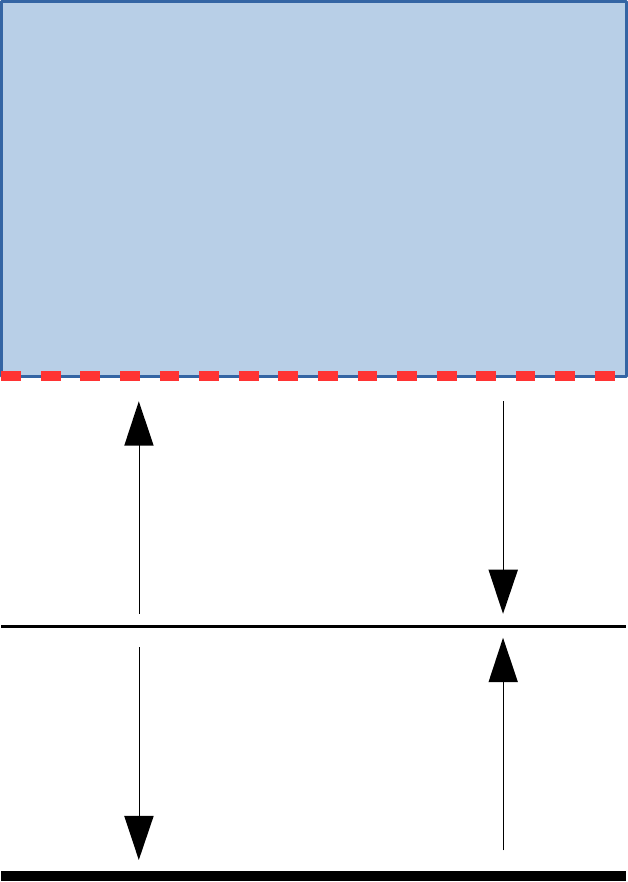}
\put(30,76){\small$\Omega_h$}
\put(50,60){\small$\partial_j \Omega_h$}
\put(75,30){\small$\Gamma_j$}
\put(75,0){\small$\Omega_l$} 
\put(60,12){\small$\Pi_{L^2(\Omega_h)}$}
\put(60,42){\small$\Pi_{L^2(\Omega_l)}$}
\put(20,12){\small$\Pi_{N(\Omega_l)}$}
\put(20,42){\small$\Pi_{N(\Omega_h)}$}
\end{overpic}
\caption{Illustration of a coupling between subdomains. $\Omega_h$ and $\Omega_l$ are the higher and lower subdomains respectively, $\Gamma_j$ is the interface between the two subdomains, $\partial_j \Omega_h$ is the portion of the boundary of $\Omega_h$ as seen from $\Gamma_j$, $\Pi_{N(\Omega_k)}$ is the projection operator from coupling variables on the mortar grid to each of the subdomains degrees of freedom ($k=h,l$), and $\Pi_{L^2(\Omega_k)}$ is the projection operator from numerical variables to coupling variables. }
\label{fig:sketch_discretization}
\end{figure}
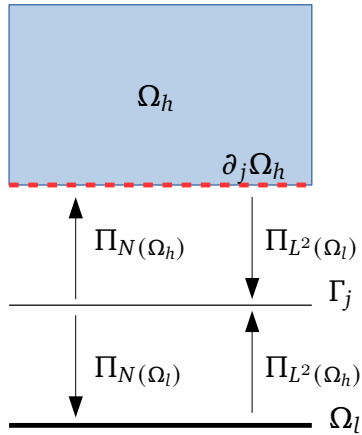

To better understand the structure of the discrete coupling, it is instructive to write out the coupled system for two subdomains $\Omega_h$ and $\Omega_l$ separated by an interface $\Gamma_j$ (see Fig. \ref{fig:sketch_discretization}). 
Let $\overline{p}_h$ and $\overline{p}_l$, be the vectors of cell-center pressures in $\Omega_h$ and $\Omega_l$ respectively, and let $\overline{\lambda}_j$ be the vector of discrete mortar fluxes in $\Gamma_j$. The discrete coupled system in absence of external sources can then be represented on the generic form
\begin{equation}
\begin{bmatrix}
A_h & 0 & G_h\Pi_{N(\Omega_h)} \\
0 & A_l & B_l\Pi_{N(\Omega_l)}+J_l\Pi_{N(\Omega_l)}T_j \\
-\Pi_{L^2(\Omega_h)} P_h & \Pi_{L^2(\Omega_l)} P_l + T_j \Pi_{L^2(\Omega_l)}R_l & D_j \\
\end{bmatrix}\begin{bmatrix}
\overline{p}_h\\
\overline{p}_l \\
\overline{\lambda}_j
\end{bmatrix}=
\begin{bmatrix}
0 \\
0 \\
0 
\end{bmatrix}.
\label{eq:matrix3p3}
\end{equation}
The first two rows of the system \eqref{eq:matrix3p3} represent the discretised differential equations in each subdomain, while the third row is the discretized Darcy's law in the direction perpendicular to the interface. Here, $A_h$ and $A_l$ are the fixed-dimensional discretizations on the subdomains, $G_h$ is the discretization of Neumann boundary conditions on $\Omega_h$, $B_l$ is the discretization of source terms in $\Omega_l$, $J_l$ is the discretization of the vector source term on $\Omega_l$, $T_j$ is the discretized $\mathfrak{K}_{\Omega\Gamma}\mathfrak{K}_{\Gamma\Gamma}^{-1}$ product on $\Gamma_j$, and $\Pi_{N(\Omega_h)}$ and $\Pi_{N(\Omega_l)}$ are the projection operators from coupling variables on the mortar grid to each of the subdomains degrees of freedom. Furthermore, $P_h$ provides a discrete representation of the pressure trace operator on $\Omega_h$, $P_l$ gives the pressure unknowns on $\Omega_l$, $R_l$ gives the reconstruction of the pressure gradient on $\Omega_l$, and $\Pi_{L^2(\Omega_k)}$ is the projection operator from numerical variables to coupling variables. Finally, $D_j$ is the discretized inverse normal permeability on $\Gamma_j$.

We conclude by making two remarks: firstly, there is no direct coupling between $\Omega_h$ and $\Omega_l$ and secondly, global boundary conditions are left out of the system.

\section{Numerical examples}
\label{sec:numerical}
We validate the semi-local model and our implementation by a suite of numerical examples. First, we consider a case with a single fault, and show how the semi-local model can capture the effects of anisotropic off-diangonal permeabilities, while the local model fails to do so.
Second, we probe the robustness of our discretization on more complex geometries in 2D and 3D.

\subsection{Comparison to the equi-dimensional model}
In this first example, we compare our reduced model to an equi-dimensional model. The aim is to highlight the enhanced modelling capabilities of our formulation with respect to the standard local formulation. With reference to this latter point, we present results of two test cases: the first one where the fault has the same off-diagonal permeability on both sides (see Fig. \ref{fig:sketch_permeabilities}.b), and a second one where different permeability structures are assigned to each side of the fault (see Fig. \ref{fig:sketch_permeabilities}.c).

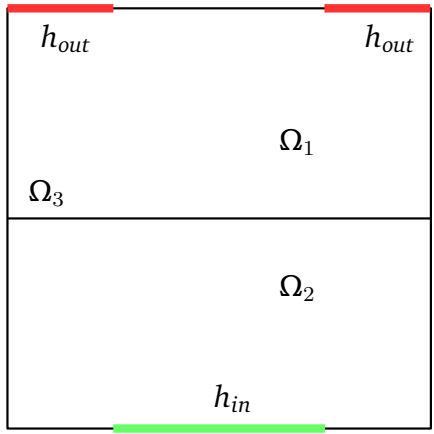
\begin{figure} []
\centering
\begin{overpic}[width=0.35\textwidth]{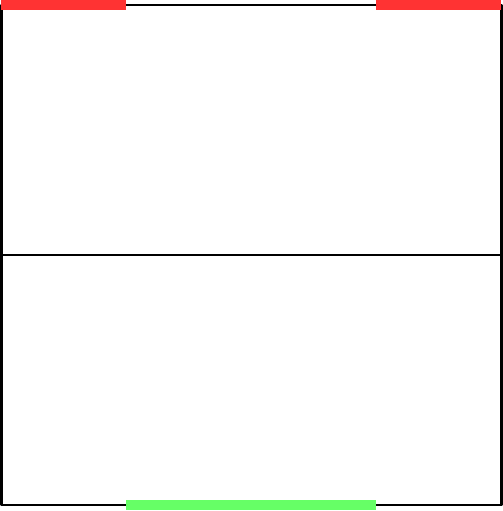}
\put(5,54){\small$\Omega_3$}
\put(63,66){\small$\Omega_1$}
\put(63,32){\small$\Omega_2$}
\put(8,90){\small$h_{out}$}
\put(83,90){\small$h_{out}$}
\put(48,6){\small$h_{in}$}
\end{overpic}
\caption{Setup of Case 1.}
\label{fig:equi_dim_homo}
\end{figure}

\subsubsection{Case 1: homogeneous permeability}
\label{sec:equi_dim_homo} % This almost sounds like Latin
We consider a 2D square domain of side $L=1~m$ cut by a horizontal fault of aperture $a=1~cm$ located in the middle of the domain. In the mixed-dimensional setting we therefore have two 2D domains $\Omega_1$ and $\Omega_2$ and one 1D fault $\Omega_3$, as illustrated in Fig. \ref{fig:equi_dim_homo}. The hydraulic conductivity is isotropic homogeneous for the 2D matrix, that is ${\bm K}_{j}=K_{j}\mathbf{I}$, with $j=1,2$, while for the fault we consider the following equi-dimensional full tensor:
\begin{equation}
{\bm K}_{3}=\begin{bmatrix}
K_{f,\parallel} & k_{f,t}\\
k_{f,t} & k_{f,\bot}
\end{bmatrix},
\end{equation} 
For simplicity, we take $K_1=K_2=K_m$. Boundary conditions consist of an applied difference in hydraulic head along the vertical direction and no-flow conditions elsewhere. In particular, the inlet pressure $h_{in}$ is specified on the portion of the bottom boundary where $0.25<x<0.75~m$, while the outlet pressure $h_{out}$ is specified on the portion of the top boundary where $x<0.25~\&~x>0.75~m$ (see Fig. \ref{fig:equi_dim_homo}). Data for the simulations are reported in Table \ref{tab:equi_dim_homo}. 
We consider as reference solution the solution obtained with an equi-dimensional model of $N=40k$ structured square cells (mesh size $dx=5~mm$), where the fault is discretized with two rows of 200 elements each. Then, for the reduced models, we consider triangular grids with approximately $N=[40, 160, 700, 3k, 11k]$ (respectively $N_f=[4, 8, 16, 32, 64]$ cells for the fault), and report the average $L^2$ error in pressure along the fault
\begin{equation}
\varepsilon_p=\dfrac{\sqrt{\sum_i \Delta_i (p_i-p_{i,eq})^2}}{\sqrt{\sum_i \Delta_i p_{i,eq}^2}},
\end{equation}
where $\Delta_i$ is the size of the fault element in the reduced model, and $p_{i,eq}$ is calculated from the equi-dimensional model as the mean value of the two fault cells at each location $x_i$:
\begin{equation}
p_{i,eq}(x_i)=\sum_{j=y_1,y_2}p_{ij},
\end{equation}
where $y_j=L/2\pm dx/2$.

\begin{figure} []
\centering
\subfloat[]{
\includegraphics[width=0.45\textwidth]{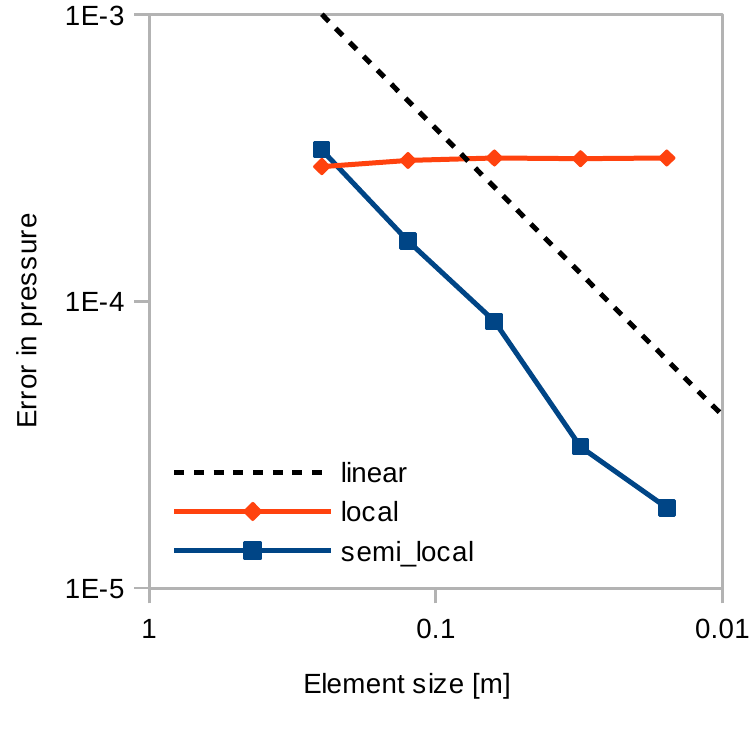}
\label{fig:convergence_equi_dim_homo}
}
\subfloat[]{
\includegraphics[width=0.45\textwidth]{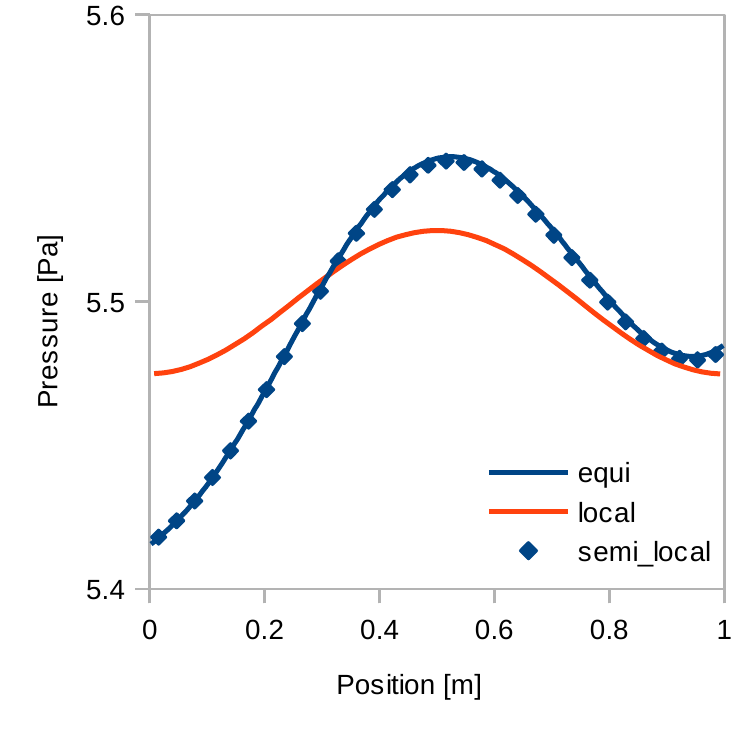}
\label{fig:pressure_along_fault}
}
\caption{Case 1: (a) convergence of the average error in pressure within the fault and (b) pressure distribution along the fault for different methods.}
\end{figure}

Convergence results are shown in Fig. \ref{fig:convergence_equi_dim_homo}. As Fig. \ref{fig:convergence_equi_dim_homo} clearly shows, our formulation presents about first-order convergence rate, while the local formulation does not converge. This is due to the strong anisotropy of the fault, which is not captured by the standard local formulation. As a result of the anisotropy of the fault, the flow will take a preferential direction towards one of the two inlets, therefore breakig the symmetry of the local formulation. This is better observed in Fig. \ref{fig:pressure_along_fault} showing the pressure distribution along the fault for the three models. As Fig. \ref{fig:pressure_along_fault} clearly shows, the semi-local and the equi-dimensional models coincide, while the local formulation exhibits an erroneous symmetric profile.

\begin{table}[]
\caption{Data for Case 1. Values of the fault hydraulic conductivity are given for the equi-dimensional model, i.e. before scaling.}
\label{tab:equi_dim_homo}
\begin{tabularx} {\textwidth} {l X r }
\hline
Parameter & Description & Value \\
\hline
{$K_m$} & Matrix hydraulic conductivity & $1~m/s$ \\
{$K_{f,\parallel}$} & Fault tangential hydraulic conductivity& $100~m/s$ \\
{$k_{f,\bot}$} & Fault normal hydraulic conductivity& $100~m/s$ \\
{$k_{f,t}$} & Fault off-diagonal hydraulic conductivity& $80~m/s$ \\
$a$ & Fault aperture & $0.01~m$\\
$L$ & Side of the square domain & $1~m$\\
$h_{in}$ & Hydraulic head at the bottom boundary & $10~m$\\
$h_{out}$ & Hydraulic head at the top boundary & $1~m$\\
\hline
\end{tabularx}
\end{table}

\subsubsection{Case 2: dual permeability}
\label{sec:equi_dim_dual_kt}
\begin{figure} []
\centering
\begin{overpic}[width=0.75\textwidth]{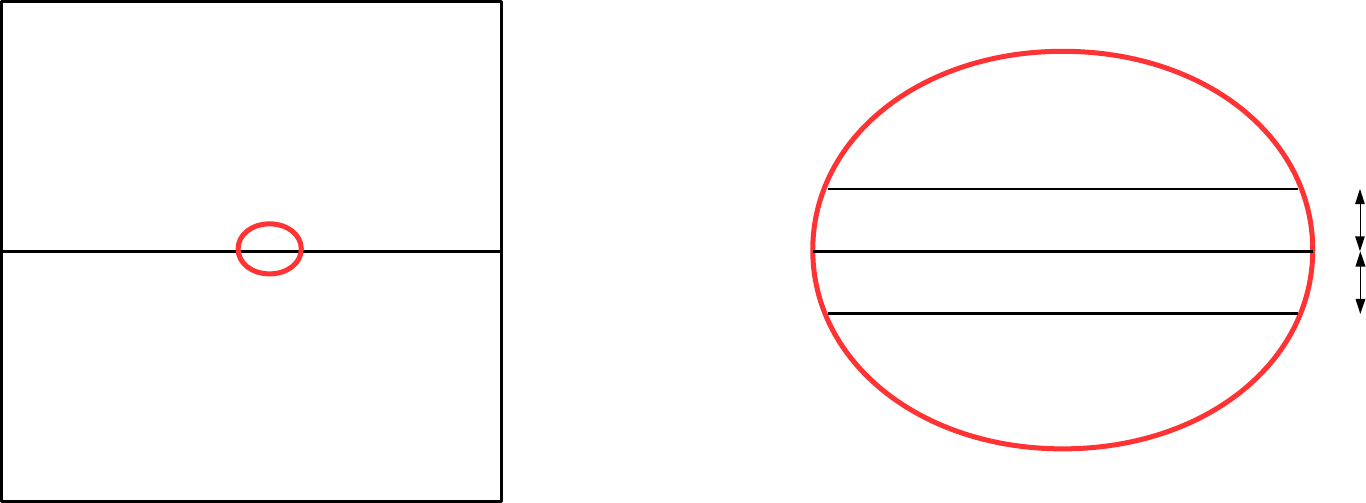}
\put(75,26){\small${K}_{1}$}
\put(75,8){\small${K}_{2}$}
\put(75,15){\small${K}_{3,2}$}
\put(75,19.5){\small${K}_{3,1}$}  
\put(102,15){$\small{1~cm}$}
\put(102,20){$\small{1~cm}$}
\end{overpic}
\caption{Setup of Case 2.}
\label{fig:sketch_dual_kt}
\end{figure}

As a further illustration of the enhanced modeling capabilities of the semi-local model, we modify the setup used in the previous section to have different permeability structures on the two sides of the fault.
This is relevant for modeling of geological faults, where the two sides of the fault may undergo different damage processes.
To that end, we divide the fault into an upper and lower part (see Fig. \ref{fig:sketch_dual_kt}) and assign different permeability structures to the two sides, that is for $j=1,2$:
\begin{equation}
{\bm K}_{3,j}=\begin{bmatrix}
K_{f,\parallel} & k_{f,j,t}\\
k_{f,j,t} & k_{f,\bot}
\end{bmatrix}.
\end{equation} 
In particular, values of $K_m$, $K_{f,\parallel}$ and $k_{f,\bot}$ are the same as those given in Table \ref{tab:equi_dim_homo}, while $k_{f,1,t}$ and $k_{f,2,t}$ take values of $50$ and $80~m/s$, respectively. The aperture of the fault is set to $a=2~cm$ and we use the same boundary conditions as in Case 1.

\begin{figure} []
\centering
\subfloat[]{
\includegraphics[width=0.45\textwidth]{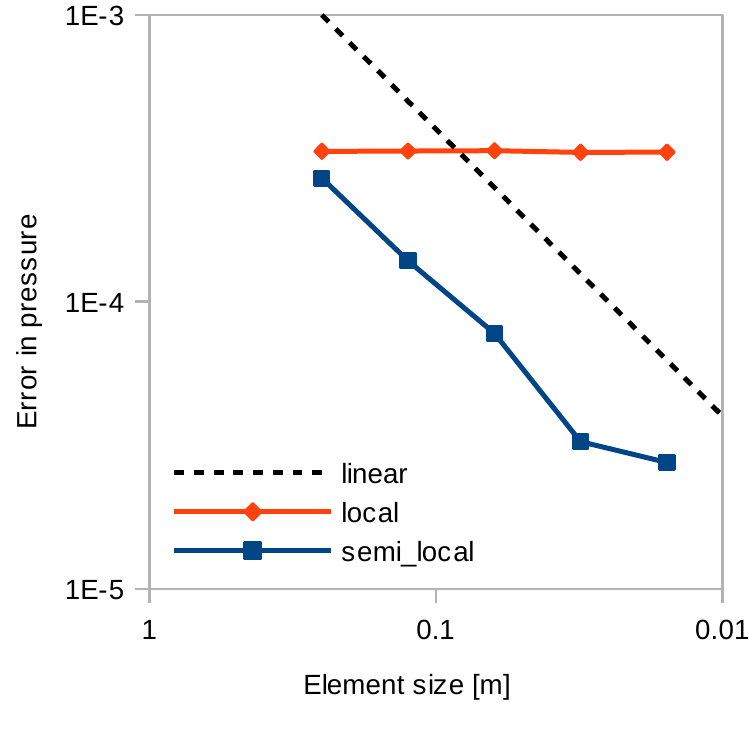}
\label{fig:convergence_equi_dim_dual_kt}
}
\subfloat[]{
\includegraphics[width=0.45\textwidth]{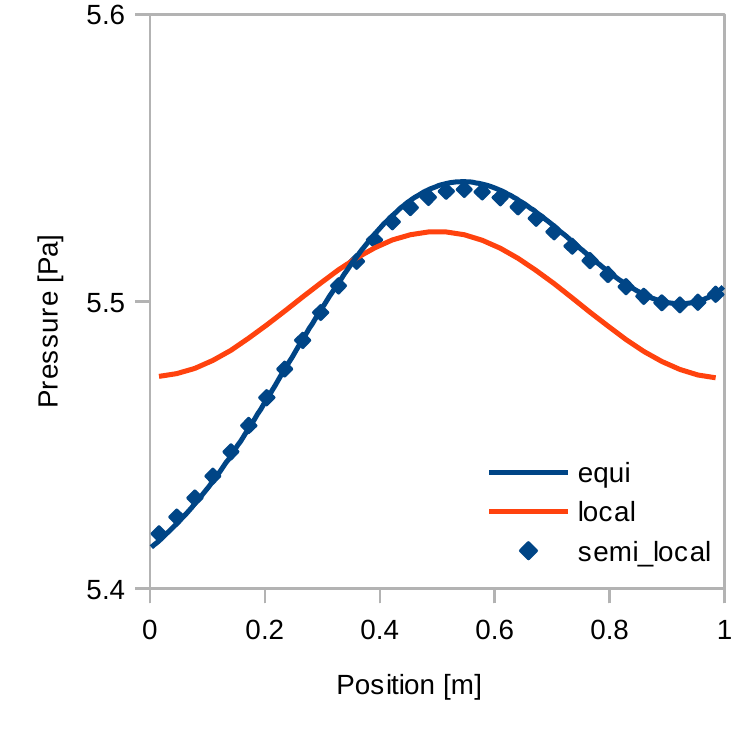}
\label{fig:pressure_along_fault_dual_kt}
}
\caption{Case 2: (a) convergence of the average error in pressure within the fault and (b) pressure distribution along the fault for different methods.}
\end{figure}

Convergence results for the local and semi-local models are shown in Figs. \ref{fig:convergence_equi_dim_dual_kt}-\ref{fig:pressure_along_fault_dual_kt}, with the reference solution again computed from an equi-dimensional model with a grid with 40k cells.
As in the previous case, the local model fails to converge, while the semi-local model exhibits first order convergence up to the last refinement step. Here, the mesh size is of the same order of the fault aperture, thus further error reduction cannot be expected due to the modeling error in the dimension reduction.

\subsection{Self-convergence}
In this section, we test the robustness of the method on more challenging fault configurations in 2D and 3D.

\subsubsection{2D case}
\label{sec:2d_intersection}
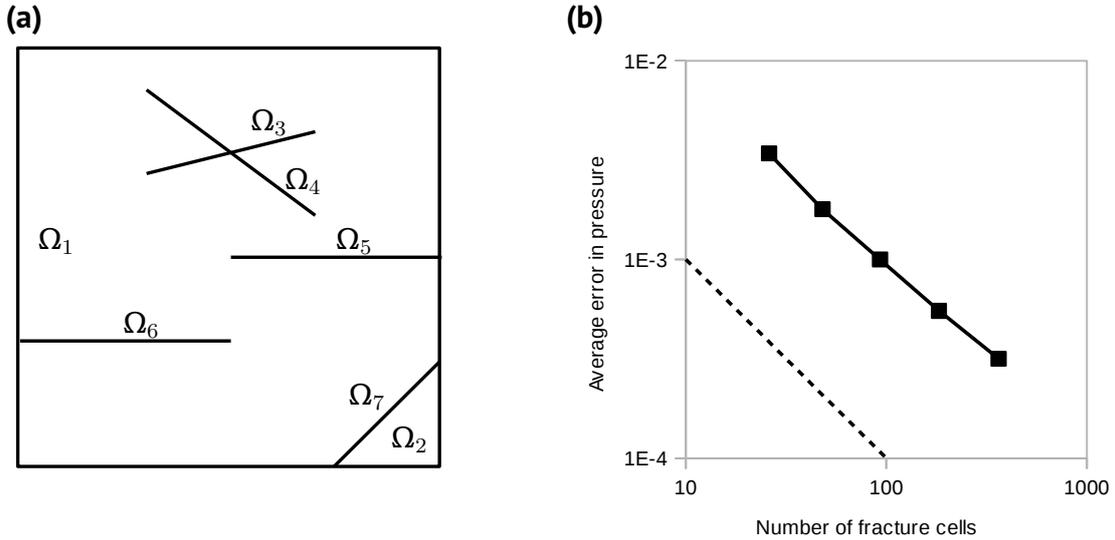
\begin{figure} []
\centering
\subfloat[]{
\label{fig:test2dintersections}
\begin{overpic}[width=0.35\textwidth]{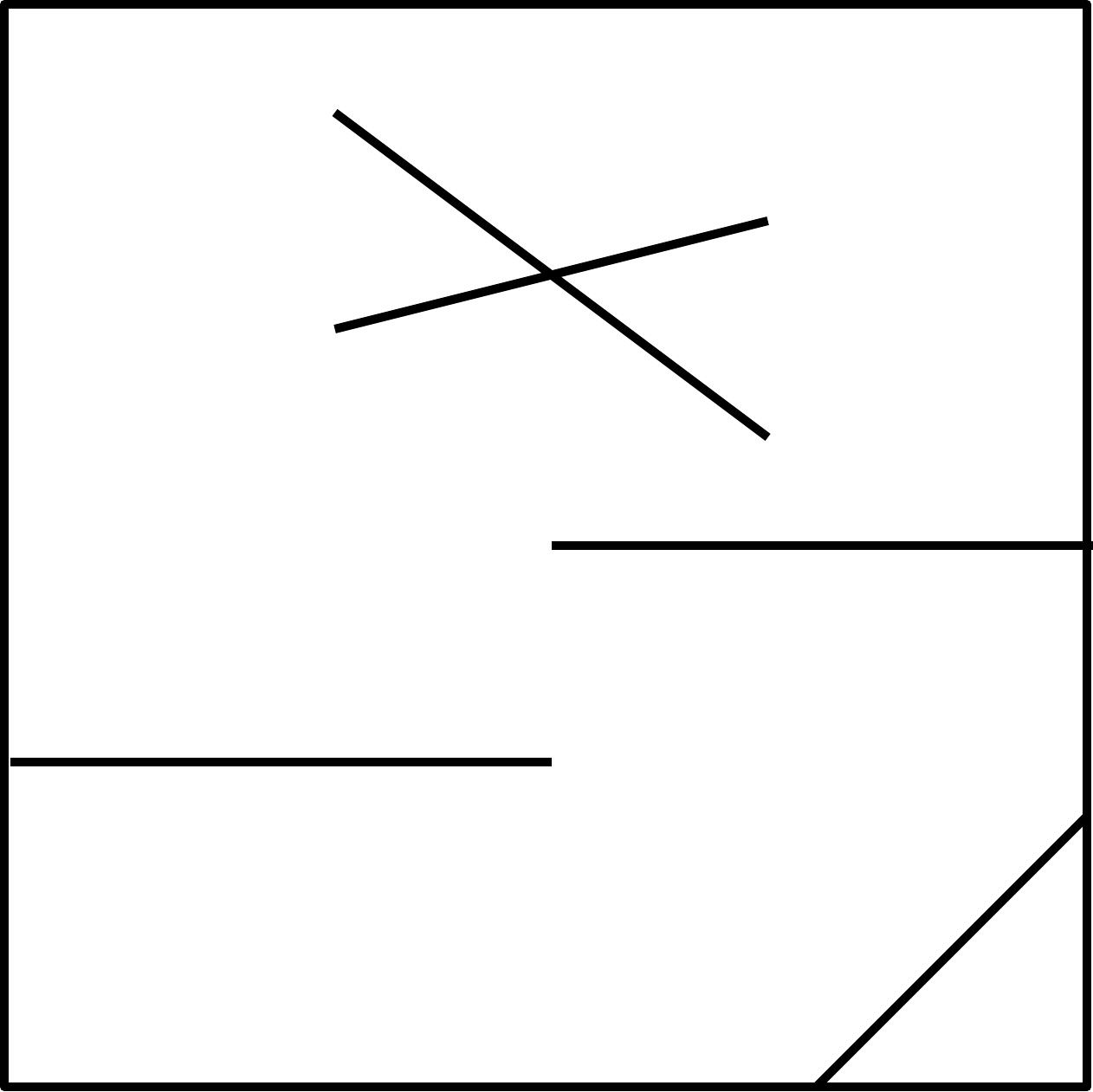}
\put(5,52){\small$\Omega_1$}
\put(88,5){\small$\Omega_2$}
\put(55,80){\small$\Omega_3$}
\put(63,66){\small$\Omega_4$}
\put(75,52){\small$\Omega_5$}
\put(25,32){\small$\Omega_6$}
\put(78,15){\small$\Omega_7$} 
\end{overpic}}
\hspace{1.5cm}
\subfloat[]{
\includegraphics[width=0.45\textwidth]{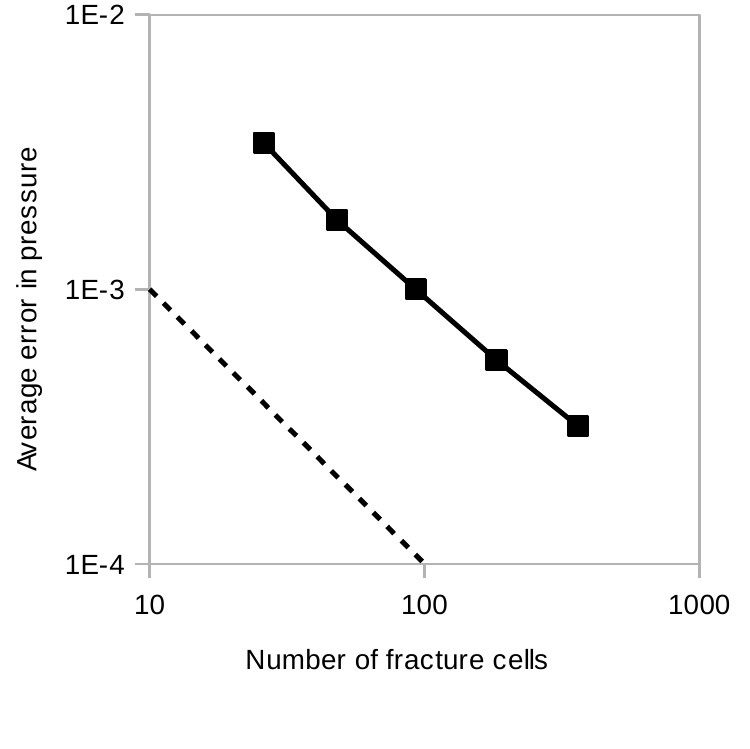}
\label{fig:convergence_2d_intersections}}
\caption{2D self-convergence test: (a) mixed-dimensional geometry and (b) convergence of the average error in pressure within the faults.}
\end{figure}

We consider the same test case as Case 1 in \citet{boon2018robust}. The domain is a unit square including a network of five faults (Fig. \ref{fig:test2dintersections}). Of these five faults, one cuts the square domain into two 2D subdomains, denoted as $\Omega_1$ and $\Omega_2$, respectively. The faults are numbered for $j=3,..,7$ and are of two kinds: $\Omega_3$ and $\Omega_4$ are conductive, that is $K_3=K_4=K_{f,1}$, while the other three are blocking, that is $K_5=K_6=K_7=K_{f,2}$. The hydraulic conductivity is isotropic homogeneous for the 2D matrix, with $K_1=K_2=K_m$, while for the faults we consider an equi-dimensional full tensor with $k_{j,t} = 0.1 K_{j,\parallel}$, for $j=3,..,7$. Boundary conditions consist of an applied difference in hydraulic head along the vertical direction and no-flow conditions elsewhere. Data for the simulations are reported in Table \ref{tab:data_2d_intersections}. 
We consider as reference solution the solution obtained with approximately $N=133k$ cells for the 2D domain and a total number of $N_f=510$ cells for the faults. Then we consider grids with approximately  $N=[300, 1k, 4k, 17k, 67k]$ (respectively $N_f=[26, 48, 93, 183, 363]$), and report the average $L^2$ error in pressure along the faults. 

The convergence results, shown in Fig. \ref{fig:convergence_2d_intersections}, indicate a rate of at least first order. The test thus confirms the performance of our method also in cases that involve faults that are intersecting and have low permeability.
Both these features are highly relevant in a geologic setting where fault may have complex geometry and reduced permeability compared to the host rock.

\begin{table}[]
\caption{Data for the 2D self-convergence test. Values of the fault hydraulic conductivity are given for the equi-dimensional model, i.e. before scaling.}
\label{tab:data_2d_intersections}
\begin{tabularx} {\textwidth} {l X r }
\hline
Parameter & Description & Value \\
\hline
{$K_m$} & Matrix hydraulic conductivity & $1~m/s$ \\
{$K_{f,1,\parallel}$} & Fault tangential hydraulic conductivity& $100~m/s$ \\
{$k_{f,1,\bot}$} & Fault normal hydraulic conductivity& $100~m/s$ \\
{$k_{f,1,t}$} & Fault off-diagonal hydraulic conductivity& $10~m/s$ \\
{$K_{f,2,\parallel}$} & Fault tangential hydraulic conductivity& $0.01~m/s$ \\
{$k_{f,2,\bot}$} & Fault normal hydraulic conductivity& $0.01~m/s$ \\
{$k_{f,2,t}$} & Fault off-diagonal hydraulic conductivity& $0.001~m/s$ \\
$a$ & Fault aperture & $0.01~m$\\
$h_{in}$ & Hydraulic head at the top boundary & $1~m$\\
$h_{out}$ & Hydraulic head at the bottom boundary & $0~m$\\
\hline
\end{tabularx}
\end{table}

\subsubsection{3D case}
\label{sec:3d_intersection}
\begin{figure} []
\centering
\subfloat[]{
\label{fig:test3dintersections}
\includegraphics[width=0.3\textwidth]{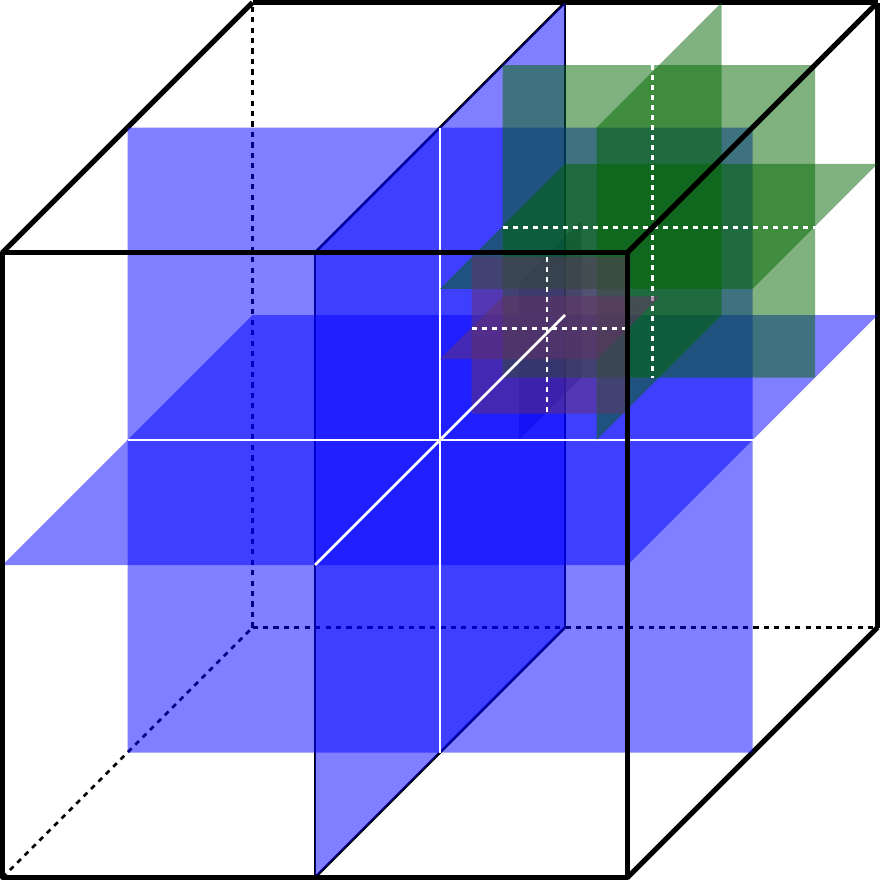}}
\hspace{1.5cm}
\subfloat[]{
\includegraphics[width=0.45\textwidth]{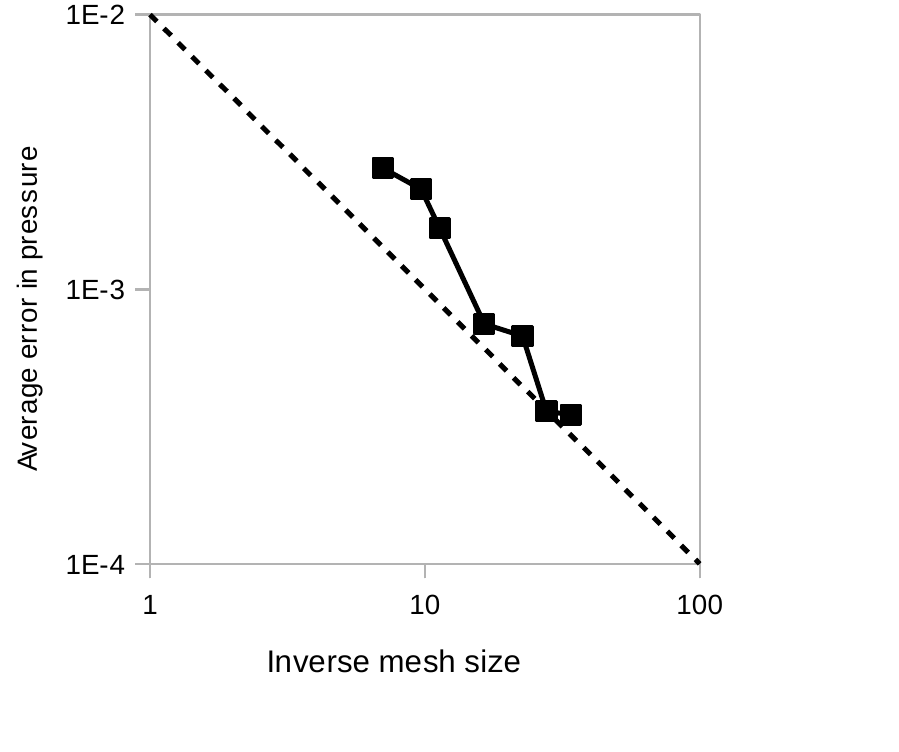}
\label{fig:convergence_3d_intersection}}
\caption{3D self-convergence test: (a) mixed-dimensional geometry and (b) convergence of the average error in pressure within the faults.}
\end{figure}

As a final verification, we consider a 3D case with multiple intersecting faults.
The setup is based on Case 2 in the benchmark study described in \cite{berre2020verification}.
The domain is a unit cube including a network of 9 faults, whose intersections divide the cubic domain into several subdomains, as illustrated in Fig. \ref{fig:test3dintersections}. 
These 3D subdomains are grouped into two regions, where we assigned 
different permeabilities $K_{m,1}$ and $K_{m,2}$, both homogeneous and isotropic (see \cite{berre2020verification} for a visualization of these two regions). For the faults we consider full tensors with tangential permeability ${\bm K}_{j,\parallel} =K_{f,\parallel} \mathbf{I}_\parallel$, normal permeability $k_{j,\bot}=k_{f,\bot}$, and off-diagonal permeability ${\bm k}_{j,t}=0.1 k_{j,t}\mathbf{i}_\parallel$. Boundary conditions consist of an imposed normal flux $q_{in}$ on the portion of the boundary where $(x,y,z)<0.25~m$ and a constant hydraulic head $h_{out}$ on the portion of the boundary where $(x,y,z)>0.875~m$. Data for the simulations are reported in Table \ref{tab:data_3d_intersections}. We consider as reference solution the solution obtained with approximately $N_3=85k$ cells for the 3D domain and a total number of $N_f=8364$ cells for all faults. Then we consider $N_3=[500, 1k, 2k, 4k, 10k, 20k, 40k]$ (respectively $N_f=[148, 282,384,814,1536,2298,3456]$) and report the average $L^2$ error in pressure along the faults. 

Convergence results are shown in Fig. \ref{fig:convergence_3d_intersection}, indicating first order convergence on average. This confirms the consistency of our implementation also for 3D problems with complex fault geometries.

\begin{table}[]
\caption{Data for the 3D self-convergence test. Values of the fault hydraulic conductivity are given for the equi-dimensional model, i.e. before scaling.}
\label{tab:data_3d_intersections}
\begin{tabularx} {\textwidth} {l X r }
\hline
Parameter & Description & Value \\
\hline
{$K_{m,1}$} & Matrix hydraulic conductivity & $1~m/s$ \\
{$K_{m,2}$} & Matrix hydraulic conductivity & $0.1~m/s$ \\
{$K_{f,\parallel}$} & Fault tangential hydraulic conductivity& $1e^4~m/s$ \\
{$k_{f,\bot}$} & Fault normal hydraulic conductivity& $1e^4~m/s$ \\
{$k_{f,t}$} & Fault off-diagonal hydraulic conductivity& $1e^3~m/s$ \\
$a$ & Fault aperture & $1e^{-4}~m$\\
$q_{in}$ & Normal flux at the inflow boundary & $-1~m/s$\\
$h_{out}$ & Hydraulic head at the outflow boundary & $1~m$\\
\hline
\end{tabularx}
\end{table}

\section{Conclusions}
\label{sec:conclusions}
We presented an improved framework to modelling and discretizing flow in generally anisotropic porous media with thin inclusions, within the context of mixed-dimensional partial differential equations. Our model considers a full permeability tensor for the inclusions, resulting in additional terms arising in our formulation as compared to existing local discretizations. We expect our model to be important for modeling of flow in faulted porous media, however the methods proposed herein can be in any case applied to models of fractures, in fact our full-permeability model naturally reduces to the existing models of fracture-matrix flow when the off-diagonal components of the inclusion permeability tensor are set to zero.
 
We provided numerical examples showing convergence of the method for both 2D and 3D faulted porous media. In particular, we provided numerical evidence that, as opposed to existing local discretizations,  our model is capable of simulating the anisotropic behaviour of the faults near damage zone. 

We remark that, in the spirit of flux-mortars coupling schemes, our formulation is independent of the discretization methods used to discretize the flow equations in the porous matrix and the faults. However, we only showed results obtained using a multi-point flux finite volume approach. Nevertheless, the formulation also applies to other discretization methods, e.g. mixed finite elements.

\textbf{Acknowledgements}\\
This work forms part of Norwegian Research Council project 250223. Data will be made public on Zenodo at the time of publication.

\bibliographystyle{apalike}
\bibliography{myref.bib} 
\end{document}